\documentclass[graybox,12pt]{svmult}

\usepackage{mathptmx}       
\usepackage{helvet}         
\usepackage{courier}        

\usepackage{makeidx}         
\usepackage{graphicx}        
\usepackage{multicol}        

\usepackage{amsmath}
\usepackage{amsfonts}
\usepackage{color}     
\usepackage{url}
\usepackage{verbatim}
\usepackage{enumerate}

\makeindex             

\usepackage{amssymb}
\usepackage{amsmath}
\usepackage{epsfig}
\usepackage{mathptmx}       
\usepackage{helvet}         
\usepackage{courier}        
\usepackage{type1cm}        

\usepackage{makeidx}         
\usepackage{graphicx}        
\usepackage{multicol}        
\usepackage[bottom]{footmisc}

\usepackage{amsbsy}
\usepackage{MnSymbol}
\usepackage{amsfonts}
\usepackage{amsopn}
\usepackage{dsfont}
\usepackage[english]{babel}
\usepackage[applemac]{inputenc}
\usepackage[colorlinks=true,breaklinks=true]{hyperref}
\usepackage[T1]{fontenc}
\usepackage{xcolor}
\usepackage{enumerate}
\usepackage{latexsym,amsmath,eufrak}
\usepackage{xcolor}
\usepackage{stackrel}
\newtheorem{teo}{Theorem}[section]
\newtheorem{prop}[teo]{Proposition}
\newtheorem{lem}[teo]{Lemma}

\newtheorem{rmk}[teo]{Remark}
\newtheorem{cor}[teo]{Corollary}
\newtheorem{defi}[teo]{Definition}
\def\i{{\infty}}

\def\i{{\infty}}

\DeclareMathOperator{\ga}{g_{\alpha,t}}
\newcommand{\limdis}{\stackrel{d}{\longrightarrow}}

\newcommand{\ID}{\mathrm{I \,D}}
\newcommand{\BO}{\mathrm{BO}}
\newcommand{\ME}{\mathrm{ME}}
\newcommand{\SD}{\mathrm{L}}

\newcommand{\STA}{\mathrm{S}}

\newcommand{\gat}{\mathbb{G}_t}

\newcommand{\ii}{{\rm i}}
\newcommand{\II}{{\rm  1~\hspace{-1.4ex}l}}
\newcommand{\nr}{\mathds{N}}
\newcommand{\zr}{\mathds{Z}}
\newcommand{\Ren}{\mathds{R}}
\newcommand{\RP}{\Ren_+}
\newcommand{\oi}{(0,\i)}

\newcommand{\Lap}{\mathcal{L}}

\newcommand{\simdis}{\stackrel{{\rm d}}{=}}

\newcommand{\CM}{\mathcal{CM}}
\newcommand{\BF}{\mathcal{BF}}

\newcommand{\CB}{\mathcal{C\!B}}

\newcommand{\TB}{\mathcal{T\!B}}

\newcommand{\LE}{\mathcal{LE}}

\newcommand\NN{\nr}

\newcommand{\qq}{\mathrm{q}}

\newcommand{\dd}{\mathrm{d}}
\newcommand{\aaa}{\mathrm{a}}
\newcommand{\bb}{\mathrm{b}}
\newcommand{\cc}{\mathrm{c}}

\newcommand{\stp}{\odot}

\DeclareMathOperator{\er}{\mathbf{E}}
\DeclareMathOperator{\pr}{\mathbf{P}}

\newcommand{\eqd}{\stackrel{d}{=}}

\newcommand\rsetminus{\mathbin{\mathpalette\rsetminusaux\relax}}
\newcommand\rsetminusaux[2]{\mspace{-4mu}
  \raisebox{\rsmraise{#1}\depth}{\rotatebox[origin=c]{-20}{$#1\smallsetminus$}}
 \mspace{-4mu}
}
\newcommand\rsmraise[1]{%
  \ifx#1\displaystyle .8\else
    \ifx#1\textstyle .8\else
      \ifx#1\scriptstyle .6\else
        .45%
      \fi
    \fi
  \fi}

\newenvironment{proofof}{\noindent {\bf Proof of }}{\ \ \ $\square$\\}
\newenvironment{proofp}{\noindent {\bf Proof. }}{\ \ \ $\square$\\}

\begin{document}
\title*{Some characterizations  of multiple selfdecomposability,  with extensions and an application to the Gamma function}
\author{
Wissem Jedidi
\thanks{Department of Statistics \& OR, King Saud University, P.O. Box 2455, Riyadh 11451, Saudi Arabia, \email{wjedidi@ksu.edu.sa}, \email{j.alromian@gmail.com}}
\thanks{Universit\'e de Tunis El Manar, Facult\'e des Sciences de Tunis, D\'epartement de Math\'ematiques, Laboratoire d'Analyse Math\'ematiques et Applications LR11ES11.  2092 - El Manar I, Tunis, Tunisia},
Zbigniew J. Jurek\thanks{University of Wroc$\l$aw, Institute of Mathematics, Pl. Grunwaldzki 2/4, 50-384 Wroc$\l$aw, Poland \email{zjjurek@math.uni.wroc.pl}},
Jumanah Al Romian$^{\star}$
}

\maketitle

\abstract{Inspirations for this paper can be traced to Urbanik \cite{urb72b,urbanik1} where convolution semigroups of multiple decomposable distributions were introduced.  In particular,  the classical  Gamma $\gat$ and  $\log \gat$, $t>0$, variables are selfdecomposable $\big($i.e. have distributions in $\SD_0(\Ren)\big)$. In fact, we show that  $\log \gat$ is twice selfdecomposable $\big($i.e. have distributions in $\SD_1(\Ren)\big)$ if, and only if, $t\ge t_1$ where $t_1$ is an explicit critical value, and this an answer to a problem raised by Akita \& Maejima  \cite{makoto}. Moreover, we provide  several new factorizations of the Gamma function and the Gamma distributions that extend many known ones in the literature. To this end, we revisit the class of multiple selfdecomposable distributions, denoted $\SD_n(\Ren)$, and propose handy tools for its characterization, mainly based on the Mellin-Euler operator and on the Hadamard  fractional integral. Finally, we give a perspective for the generalization of the class $\SD_n(\Ren)$, based on linear operators or on stochastic  integral  representations.}

\smallskip

\begin{center}
{\it In memoriam of Kazimierz Urbanik (February 5, 1930 -- May 29, 2005)}
\end{center}

\section{Introduction}
The distribution $\mu$  of a real-valued random variable $X$, is said to be  {\it infinitely divisible}, and we denote $\mu \in \ID(\Ren)$ or $X\sim \ID(\Ren)$, if
$$\mbox{for each natural $n\ge 2$ there exits $\mu_n$ an $n$-th fold convolution: $\;\mu_n^{\ast n}=\mu.$}$$
In other terms, $\mu_n$  is the distribution of independent and identically distributed random variables $X_{1,n},\, X_{2,n}, \, \ldots , X_{n,n}$ and
$X \eqd  X_{1,n} + X_{2,n} + \ldots + X_{n,n}.$  
It is known  that every r.v. $X\sim \ID(\Ren)$ is embedded into a {\it L\'evy process} $(Z_t)_{t \geq0}$, i.e. a process with independent and stationary increments, such that $Z_0=0$ and  $X \eqd  Z_1$, in order that the so-called L\'evy-Khintchine formula holds:
\begin{equation}\label{zt}
X\sim \ID(\Ren) \Longleftrightarrow \er[e^{\ii  \,u\,  X}]^t = \mathbb{E}[e^{\ii \,u\,Z_t}]=e^{t\, \Phi (u)}, \quad t>0,  \;\; u \in \Ren,
\end{equation}
where $\Phi$ is given by the following expression: for a fixed truncation function $h$ (i.e. a bounded function such that $\lim_{x\to 0} (h(x)-x)/x^2$ exists, for example $h(x)=x \II_{|x|\le 1}$ or $x/(1+x^2)$, we have
\begin{equation}\label{LK}
\Phi(u)= \ii\, \aaa \, u - \bb\, u^2 + \int_{\Ren\rsetminus \{0\}}\big(e^{\ii \,u\,x}-1-\ii \, u \, h(x)\big) \pi(dx),
\end{equation}
where $\aaa \in \Ren, \, \bb\geq0 $ are respectively called {\it drift  term} and {\it Brownian  coefficient}. Clearly, if one changes the truncation function to a new one $g$, the new drift term becomes $\aaa_g:=\aaa + \int_{\Ren\rsetminus \{0\}} \big(g(x)-h(x)\big) \pi(dx)$. The measure $\Pi$ is called the {\it L\'evy measure} and satisfies the integrability condition
\begin{equation}\label{x21}
\int_{{\Ren\rsetminus \{0\}}}(x^2 \wedge 1) \pi(dx)< \infty.
\end{equation}
\begin{itemize}
\item {\it Spectrally negative} L\'evy processes on the real line have a particular interest. They  correspond  to  L\'evy processes $(Z_t)_{t \geq0}$ with no positive jumps, viz. such that the associated L\'evy measure $\Pi$ in \eqref{LK} gives no mass to $(0, \infty)$. Defining
\begin{equation}\label{image}
 \mbox{$\Pi$  be the image of $\pi$ by the reflection $x\mapsto -x,\quad$ and $\quad \chi(x):=-h(-x), \quad \mbox{$h$ given by \eqref{LK}}$},
\end{equation}
it becomes more handy to use the so-called {\it L\'evy-Laplace exponent} $\Psi$ instead of  the {\it L\'evy-Fourier} exponent $\Phi$:
\begin{equation}\label{ird}
\Psi(\lambda):= \lim_{t\to  0} \frac{\log \er[e^{\lambda\, Z_t }]}{t}= \Phi(-\ii \lambda)= \aaa \lambda \,  + \bb\, \lambda^2 + \int_{(0,\infty)}\big(e^{-\lambda x}-1 +\lambda \, \chi(x)\big) \Pi(dx), \quad \lambda \geq 0,
\end{equation}
where $\Pi$ satisfies the integrability condition \eqref{x21}. The class $\LE$ of {\it Laplace exponents} is defined as the set of functions $\Psi$ of the form \eqref{ird} with $\chi(x)=x$, i.e.
\begin{equation}\label{lev}
\LE:=\big\{\Psi(\lambda)=   \aaa \lambda + \bb \lambda^{2}   + \int_{(0,\infty)} \big(e^{-\lambda x} - 1
+ \lambda x \big)\,\Pi (dx), \quad \lambda \geq 0\big\},
\end{equation}
where $\aaa\in \mathds{R}$, $\bb \geq 0$, and the L\'evy measure $\Pi$ satisfies integrability condition
\begin{equation}\label{integg}
\int_{(0,\infty)}  (x^2\wedge x)\Pi(dx)<\infty.
\end{equation}
We define the subclass $\overline{\ID}_-(\Ren)$ and $\ID_-(\Ren)$ of $\ID(\Ren)$ as follows:
\begin{equation}\label{fold1}
X\sim \overline{\ID}_-(\Ren) \; \big(\mbox{respectively} \;\ID_-(\Ren)\big) \quad \mbox{if}\quad  \er[e^{\lambda X}]^t =\er[e^{\lambda Z_t}]=e^{t\Psi(\lambda)}, \;\; t>0 \quad \mbox{and $\Psi$ is as in $\eqref{ird}\; \,\big($respectively $\Psi \in \LE\big)$.}
\end{equation}
For L\'evy processes, the references are numerous, we suggest the books of Bertoin \cite{bert}, Kyprianou \cite{kyp} or Sato \cite{sato}. The following  proposition explains why $\overline{\ID}_-(\mathds{\Ren})$ is the vague the closure of $\ID_-(\Ren)$:
\begin{prop} We have the equivalences:
\begin{enumerate}[1)]
\item $X\sim \overline{\ID}_-(\mathds{\Ren}) \Longleftrightarrow  X \eqd \lim_{n} X_n,\;$ with $X_n\sim \ID_-(\Ren)$.
\item $X\sim \ID(\Ren) \Longleftrightarrow $ there exist two sequences of  independent r.v's $X_n', \,X_n''\sim \ID_-(\Ren)$ such that $X \eqd  \lim_{n} (X_n'- X_n'')$.
\end{enumerate}
\label{vague}\end{prop}
From the latter, we see that $X\sim \ID(\Ren)$ if, and only if $X=X'-X''$, where $X'$ and $X''$ are independent and $X', \; X''\sim \overline{\ID}_-(\mathds{\Ren})$.

\item The class $\BF$ of {\it Bernstein functions} is defined by
\begin{equation}\label{lev1}
\BF :=\big\{\phi(\lambda)= \dd \lambda   + \int_{(0,\infty)} \big(1-e^{-\lambda x} \big)\,\Pi (dx), \quad \lambda \geq 0\big\},
\end{equation}
where $\dd\geq 0$ and  the L\'evy measure $\Pi$ satisfies the integrability condition
\begin{equation}\label{x11}
\int_{(0,\infty)} (x  \wedge 1)\Pi (dx) <\infty.
\end{equation}
We emphasize that, traditionally, a function of the form $\phi+\qq$ is called Bernstein function, when $\phi \in\BF$ and the addition of so-called {\it killing term} $\qq\geq 0$, transforms probability measures into sub-probability measures. In order to avoid discussions on sub-probability measures, we  only  consider the class $\BF$. Observe that a differentiable function $\Psi$ on $(0,\infty)$  belongs to $\in \LE$ if, and only if, $\er[Z_1]=\Psi'(0)$ is finite and  $\Psi'-\Psi'(0)\in \BF$.  Also observe that  and if that $\phi \in \BF$ is equivalent to the fact that $\phi$ is nonnegative and $\phi'$ belongs to the class of {\it completely monotone} functions and this yields  $-\phi \in \LE$. The class $\CM$ of completely monotone functions corresponds to those infinitely differentiable functions $f:\oi \to \oi$,  such that
$$(-1)^n f^{(n)} \geq 0,\quad  \text{for all  $n=0,1,2\cdots\;$}.$$
Bernstein characterized the class $\CM$  by
\begin{equation}\label{cmr}
\CM:=\big\{f(\lambda)=\int_{[0,\infty)}e^{-\lambda x} \nu (dx), \quad \lambda >0, \quad \mbox{where $\nu$ is some Radon measure on $[0,\infty)$}\big\}.
\end{equation}
As for L\'evy-Laplace exponents, the class $\BF$ is one-to-one with the set
of infinitely divisible distributions supported by $\Ren_+$:
$$\mbox{we denote $X\sim \ID(\Ren_+), \; \;$  if $\; \;X\geq 0$ and $X\sim \ID(\Ren)$,}$$
which is equivalent to $X\geq 0$   and $-X\sim \ID_-(\Ren)$. In this case, the embedding L\'evy process $(Z_t)_{t \geq0}$  in \eqref{zt} is a {\it subordinator} (i.e. a L\'evy process with increasing paths) and the  L\'evy-Khintchine formula takes the form
\begin{equation}\label{fold}
X\sim \ID(\Ren_+) \Longleftrightarrow  \er[e^{-\lambda X}]^t =\er[e^{-\lambda Z_t}]= e^{-t\phi(\lambda)}, \;\; t\geq 0, \;\; \lambda  \geq 0, \quad \mbox{and}\; \; \phi \in \BF.
\end{equation}
 \end{itemize}
The class $\SD_0(\Ren)$  of {\it selfdecomposable} distributions has the following description, among others:
\begin{equation}\label{hhhh}
\mbox{we denote}\;\mu\in \SD_0(\Ren), \quad \mbox{if, for all $c\in (0,1)$, there exists a measure $\mu_c$ s.t.}  \;\; \mu=T_c\mu \ast \nu_c.
\end{equation}
where, for Borel sets $B, \;T_c\mu(B)= \mu(c^{-1}B)$. The latter has several equivalent formulations that can be found in \cite{JV1983} and in the book of Sato \cite{sato}:
\begin{eqnarray}
X\sim  \SD_0(\Ren) &\Longleftrightarrow& X\simdis c\, X +  V_c, \quad \mbox{where $V_c$ is independent of $X$ (necessarily  $V_c \sim \ID(\Ren)$),} \label{xcy}\\
&\Longleftrightarrow& X\sim \ID(\Ren),\;\; \mbox{and for each $0<c<1$, and each Borel set $B\subset \Ren$, the L\'evy measure satisfies}, \;\; \Pi(B)-\Pi(c^{-1}B)\ge 0, \nonumber\\
&\Longleftrightarrow& \Pi(dx)=\frac{k(x)}{x} dx, \quad \mbox{where the function $k$ is non-increasing on $(-\i, 0)$ and on $(0,\i)$}.\label{sdp}
\end{eqnarray}
Observe that
$$X\sim  \SD_0(\Ren) \;\; \mbox{and} \;\; X \geq 0 \Longleftrightarrow \mbox{the $k$-function in  \eqref{sdp} is such that $\;k(x)=0$, for $x<0$}.$$

\medskip

In 1973, Urbanik \cite{urb72b,urbanik1}, introduced a decreasing sequence $\SD_n(\Ren),\; n=1,\ldots$, of limiting distributions (for some specified triangular arrays) contained in the class $\SD_0(\Ren)$, that is,
\begin{equation}
\SD_\infty(\Ren) =\bigcap_{k\geqslant 0} \SD_k(\Ren) \subset \ldots \subset\SD_n(\Ren)  \subset\ldots \subset\SD_1(\Ren) \subset  \SD_0(\Ren).
\end{equation}
Probability  measures in $\SD_n(\Ren)$ are called \emph{$n$-times selfdecomposable measures.}
Urbanik characterized  multiple selfdecomposable distributions  by their characteristic functions in \cite{urbanik} and \cite[Theorem 1 and 2]{urbanik1}. Then Urbanik gave in \cite{urb72b} the results without detailed  proofs. His proofs in \cite{urbanik,urb,urb1972} used the Choquet-Krein-Milman theorem on extreme points in compact convex sets; as a reference for this  cf. Phelps \cite{Ph}. A probabilistic proof, using random integral representations, is in \cite{JU2}. Urbanik's  description of the classes $L_n(\Ren)$, in terms of the convolution factorization is given in \cite[Proposition 1]{urbanik1}: with the convention $\SD_{-1}(\Ren):=\ID(\Ren)$,
\begin{equation}
\mu\in \SD_n(\Ren), \quad  \mbox{if, for all $c\in (0,1)$, there exists $\mu_c\in \SD_{n-1}(\Ren)\;\;$ s.t.}  \;\; \mu=T_c\mu \ast \nu_c,
\end{equation}
For more convenience, we denote
$$\SD_n(\mathds{\Ren_+}):=\big\{\mu \in \ID(\Ren_+);\;\mbox{\it such that}\;\mu\in \SD_n(\Ren)\big\}.$$
Note that $\ID(\mathds{A})$ and $\SD_n(\mathds{A}), \;\mathds{A}= \Ren_+, \;\Ren$, are closed (in the weak convergence topology) convolutions
semigroups.   If $\mu\in \SD_n(\Ren) \setminus \SD_{n+1}(\Ren)$, then we say that $\mu$ is \emph{exactly $n$-times selfdecomposable.}  For operator-selfdecomposability problems, we refer the book of Jurek \& Mason \cite{JU} and Meerschaert \& Scheffler  \cite{Meer}, for distributional properties of selfdecomposable distributions we recommend the book of Sato \cite{sato}. For analytic properties related to infinitely divisible distributions on the half-real line, we recommend the book of Schilling, Song \& Vondra\v{c}ek \cite{SSV} and also  Steutel \& van Harn \cite{steutel}. Several other proofs for the characterization of   $\SD_n(\Ren_+)$ are $\SD_n(\Ren)$ are available and the references are multiple, cf. Berg \& Forst \cite{berg}, Jurek \cite{JU2, JU0, JU, JUS}, Van Thu \cite{van} and also Steutel \& van Harn \cite{steutel} and the references therein.\\

In Section \ref{newg}, we provide several new  decomposability properties for the  Gamma function, the Gamma distributions, and the positive stable distribution. For this purpose, we consider, for  $\alpha \in (0,1)\cup(1,\infty)$  and  $ t> 0$, the following functions:
\begin{equation}\label{nota}
 G_{\alpha }(\lambda)=\frac{\Gamma(\lambda)^\alpha }{\Gamma(\alpha \lambda)}, \; \; \lambda>0 \quad \mbox{and}\quad
G_{\alpha ,t}(\lambda)= \Gamma(t)^{1-\alpha} \; \alpha^{\alpha \lambda} \, \,\frac{\Gamma(\lambda +t)^\alpha }{\Gamma(\alpha \lambda + t)},\;\; \lambda \geq 0,
\end{equation}
and notice the following relation
\begin{equation}\label{nice}
G_{\alpha ,t}(\lambda)= \frac{1}{G_{\frac{1}{\alpha} ,t}(\alpha \lambda)^{\alpha}}.
\end{equation}

\begin{enumerate}
\item Theorem \ref{gamXY} retrieves Berg, \c{C}etinkayab \& Karp's result \cite[Theorem 3.13]{bergkarp} and gives a stochastic interpretation. The latter and Theorem  \ref{gself}  provide  the exact range of the parameters $\alpha$ and $t$  for which
$$G_{\alpha ,t} \;\;  \left(\mbox{respecively}\;\; 1/G_{\alpha ,t} \right)
 \;\;\mbox{is the Laplace transforms of  a r.v  $X_{\alpha ,t}$ (respecively $Y_{\alpha ,t}$) with distribution in  $\ID(\Ren_+)$ and $\SD_0(\Ren_+)$}.$$
\item Corollary \ref{pidg}  proves that if $\alpha _1,\alpha _2,\ldots ,\alpha _n \in(0,1)$,  $\sum_{k=1}^n\alpha_k=1$, and if $\gat,  \mathbb{G}_{1,t},\,\ldots ,\,\mathbb{G}_{n,t}$, are independent random variables with Gamma distribution with shape parameter $t>0$, then
$$\gat \simdis d({\underline{\alpha}})\;\, \mathbb{G}_{1,t}^{\alpha _1}\,\ldots \,\mathbb{G}_{n,t}^{\alpha _n} \,\;e^{-X_{\underline{\alpha},t}} \Longleftrightarrow t \geq 1/2 \quad \mbox{and}
\quad d({\underline{\alpha}})\;\, \mathbb{G}_{1,t}^{\alpha _1}\,\ldots \, \mathbb{G}_{n,t}^{\alpha _n} \simdis \gat  \,\;  e^{- Y_{\underline{\alpha},t}}  \Longleftrightarrow t < 1/2,$$
where $d(\underline{\alpha}):= \prod_{k=1}^n\alpha_k^{\alpha_k}$ and all r.v.'s are assumed to be independent in each side of the identities and $X_{\underline{\alpha},t}, \;Y_{\underline{\alpha},t} \sim\ID(\Ren_+)$. Analytically, the latter is transcripted in terms of the Gamma function with  the following remarkable fact:
  \begin{equation}\label{rab}
  \lambda \mapsto 
\left(\frac{\Gamma(\lambda +t)}{d(\underline{\alpha})^\lambda\prod_{k=1}^n \Gamma(\alpha_k \lambda +t)}\right)^r \; \mbox{is completely monotone},
  \end{equation}
if, and only if $t\geq \frac{1}{2}$ and $r>0$ or  $t<\frac{1}{2}$ and $r<0$. This results generalizes several known results in the literature, c.f Alzer \& Berg \cite{Alzer}, Bertoin \& Yor \cite{BY}, Li \& Chen \cite{lichen},  Mehrez \cite{tun},  Pestana,   Shanbhag \& Sreehari \cite{PSS}, for instance. We also explore whenever $X_{\underline{\alpha},t},\; Y_{\underline{\alpha},t}\sim \SD_0(\Ren_+)$ or $\BO$, ($\BO$ being the Bondesson class of distributions, cf. \cite[Definition 9.1]{SSV}).
\item We improve Akita \& Maejima's \cite{makoto} results who proved that $\log \gat$ is twice selfdecomposable for any $t\geq 1/2$ and that there exists a universal constant $t_1\in (0,1/2)$  such that the last property could be extended for $t \in (t_1, +\infty)$. Using  tools developed in Section \ref{simple}, we were able in Proposition \ref{maj} to show that the value of $t_1$ is the maximum of an explicit elementary function, and is approximately 0.151649938034.
\end{enumerate}

Section \ref{Mellic}  prepares for the results of Section \ref{simple} and Lemma \ref{chark} gives the full characterization of the class
\begin{equation} \label{Mn}
\mathcal{M}_n:= \big\{k, \; \mbox{s.t. $\;x\mapsto k(e^x)\;$ is $n$-monotone  on $\Ren$}\big\},
\end{equation}
cf. Williamson \cite{wil} for $n$-monotone functions. Using the so-called {\it Euler-Mellin differential operator} $\Theta =x\;d/dx$ and the so-called Hadamard
fractional integral, we provide the converse of \cite[Proposition 1.16]{SSV} and, the analog of \cite[Theorem 4.11]{SSV}  for the class $\mathcal{M}_n$. \\

In Section \ref{simple}, mainly in Corollary \ref{char}, we shall provide a simple proof for the characterization of the classes $\SD_n(\Ren)$ based on the class $\mathcal{M}_n$ given in \eqref{Mn}.  Jurek had provided the same characterization \cite[Theorem 6.2 and Theorem 7.1]{{JU1}} for infinitely divisible measures on Banach spaces.
Our approach  may have the merit to exhibit the tight link between the Mellin-Euler differential $\Theta$ operator and the classes $\SD_n(\Ren)$ and suggests that a mechanism producing new classes of infinitely divisible distribution via linear operators $\Omega$ other than $\Theta$ could be implemented as follows: assume $\Omega$ is a linear operator on the space of functions $f:[0,\infty)\to [0,\infty)$ which commute with the dilations, i.e. there exists
$$\Omega(x\mapsto f(ux))= \Omega(f)(ux) , \quad \mbox{for all} \; x,\;u>0.$$
For instance, if $\Omega=\Theta$,  then $\Omega(\BF)\cap \BF$  is one-to-one with $\SD_n(\Ren_+)$ and if  $\Omega=I-\Theta$, then $\Omega(\BF)\cap \BF$ is not void.  Now, let $\phi$ be a  Bernstein function of the form
$$\phi(\lambda)= \dd\; \lambda   + \int_{(0,\infty)}  (1-e^{-\lambda x})\;\;\frac{k(x)}{x} dx =\dd \lambda   + \int_{(0,\infty)} (1-e^{-x})\,\tilde{k}(\lambda/x)\, \frac{dx}{x},\quad \lambda \geq 0.$$
where $\dd\geq 0$ and assume that $\tilde{k}(x):= k(1/x)$  is ``good  enough" so that we can swap $\Omega$ and the integral. If  $\Omega^n, \; n=1,2,\ldots$, is the $n$-th iterate of $\Omega$ and  $\omega:=\Omega(Identity)(1)$, then,
$$\Omega^n(\phi)(\lambda)=   \omega\; \dd \; \lambda   + \int_{(0,\infty)}  (1-e^{-\lambda x} )\; k_n(x)\; \frac{dx}{x}, \quad\mbox{where}\quad k_n(x):=\Omega^n (\tilde{k})(1/x).$$
Thus, it is simply seen that $\Omega^n(\phi)$ remains a Bernstein function if $k_n$ complies with the non-negativity and integrability conditions imposed on L\'evy measures. With this mechanism, one is able to build new subclasses of  $\ID(\Ren_+)$, eventually nested.  When $\Omega=\Theta$, we obtain in Corollary \ref{char} the exact shape of the function $k_n$ (for Bernstein or L\'evy-Laplace exponents) of $X \sim  \ID(\Ren_+)$ or $\ID_-(\Ren)$ or  $\overline{\ID}_-(\Ren)$ for which it holds that
$$X \sim  \SD_n(\Ren_+) \; \; \mbox{or}  \; \; \SD_n(\Ren), \quad n=1,\, 2,\ldots, \infty.$$

Section \ref{verva} improves the last approach with a method based on integral stochastic  representations like  \eqref{verv}.  All the proofs were postponed to Section \ref{proo} and Section \ref{acco} gives some background on the classes $\ID(\Ren)$ and $\SD_0(\Ren)$.
\section{New  decomposability properties for the  Gamma function, the Gamma distributions and the positive stable distribution}\label{newg}
This section provides new multiple selfdecomposability properties involving the  Gamma function and the Gamma and positive stable distributions. \\

Recall the Digamma function is defined by $\psi(t)= \Gamma'(t)/\Gamma(t), \; t>0$ and is given by formula 5  p.903   \cite{grad}:
\begin{equation*}\label{psi}
\psi(t)=  -\gamma + \,\int_0^\infty  \frac{e^{- x}- e^{-t x}}{1- e^{-x}}   dx, \quad t>0,
\end{equation*}
where $\gamma$ is the Euler-Mascheroni constant.
Since $\log \Gamma(\lambda)=\int_1^\lambda \Psi(t) dt,\; \lambda>0$, we recover the following representations,
$$\Gamma(\lambda) = \exp\left\{- \gamma (\lambda -1) +\int_0^\infty \!\!( e^{-\lambda u} -e^{-u} - (\lambda-1) u e^{- u} )\,\frac{du}{u (1-e^{-u})} \right\},\quad  \lambda>0$$
and then,
\begin{equation}
\frac{\Gamma(\lambda +t)}{\Gamma(t)}=\exp\left\{\psi(t)\lambda+\int_0^\infty \!\!(e^{-\lambda u}-1 +\lambda u)\,\frac{ e^{-tu}}{u (1-e^{-u})}\;du
\right\} , \quad \lambda\geq 0,\;t>0.\label{logammat}
\end{equation}
The $q$-Gamma function $\Gamma_q(x), \; x>0$,  defined by
\begin{equation}\label{qgam}
\Gamma_q(x) =  \left\{ \begin{array}{lll}
&(1-q)^{1-x} \; \prod_{j=0}^\infty  \frac{1-q^{j+1}}{1-q^{j+x}},& \quad \mbox{if} \;0<q<1,\\
&\mbox{and}& \\
&(q-1)^{1-x} \;q^{x(x-1)/2}\; \prod_{j=0}^\infty  \frac{1-q^{-(j+1)}}{1-q^{-(j+x)}},&  \quad \mbox{if} \;q>1,
\end{array}\right.
\end{equation}
enjoys the basic property: $\lim_{q\to 1-} \Gamma_q(x) = \lim_{q\to 1+} \Gamma_q(x) = \Gamma(x)$. \\

From now on, we denote by $\gat$ a random variable with the {\it standard  Gamma distribution}  with shape parameter  $t >0$, which  has the density function, Laplace and Mellin transforms respectively given by
\begin{equation}  \label{logammatt}
f_{\gat}(x)= \frac{x^{t-1}}{\Gamma(t)}\, e^{-x}, \quad \er[e^{-\lambda\gat}] =\frac{1}{(1+\lambda)^t} \quad\mbox{and}\quad \er[\gat^\lambda] ={\Gamma(t+\lambda) \over
\Gamma(t)},\quad  x>0, \;\;\lambda>-t.
\end{equation}
The function $ \lambda \mapsto \lambda^\alpha, \, 0<\alpha <1,$ is a generic  example of a Bernstein functions. It is not difficult to derive the representation
\begin{equation}\label{stabphi}
\lambda^\alpha = \frac{\alpha}{\Gamma (1-\alpha)}\int_0^{\infty} (1-e^{-\lambda x}) \frac{dx}{ x^{\alpha+1}} ,\quad  \lambda \geq 0,
\end{equation}
and this function is associated to the so-called standard {\it positive stable} r.v.  $\mathbb{S}_\alpha$ with stability parameter  $\alpha \in (0,1)$ via:
\begin{equation}\label{stap}
\er[e^{-\lambda \mathbb{S}_\alpha}]= e^{-\lambda^\alpha},\quad \lambda \geq 0, \qquad \mbox{and}\qquad \er[(\mathbb{S}_\alpha)^{-\lambda}]=
\frac{\Gamma(1+\frac{\lambda}{\alpha}) }{\Gamma(1+\lambda)}, \quad \lambda >-\alpha.
\end{equation}
The p.d.f. of $\mathbb{S}_\alpha$  is not explicit, except for $\alpha =1/2$  where $S_{1/2} \simdis 1/(4 \;\mathbb{G}_{1/2})$, cf. the monograph of Zolotarev \cite{zolo} for more account.
\subsection{Elementary  decomposability properties for Gamma and stable distributions}
Using  the   representations \eqref{logammat} and \eqref{logammatt}, we deduce that the functions, given for $\lambda \geq 0$, by
\begin{eqnarray}
\phi_{\gat}(\lambda)&=&-\log \er[e^{-\lambda \gat}]=  \int_0^\infty \big(1-e^{-\lambda x} \big)\,\frac{k_1(x)}{x} \,dx,
\quad k_1(x)=e^{-x},\nonumber \\
\Psi_{\gat}(\lambda)&=& \log \er[(\gat)^\lambda]=  \int_0^\infty \big(e^{-\lambda x} - 1+\lambda x \big)\,\frac{k_2(x)}{x} \,dx,
\quad k_2(x)= \frac{e^{-tx}}{1-e^{-x}},\label{logamma1}
\end{eqnarray}
are respectively in $\BF$ and $\mathcal{LE}$, and that $k_1, \, k_2$  are non increasing (cf. also \cite[p. 98]{JU10} for $k_2$). By \eqref{sdp}, it is immediate  that $\gat \sim \SD_0(\Ren_+)$ and $\log \gat \sim \overline{\ID}_-(\Ren)\cap \SD_0(\Ren)$; cf. Jurek \cite[p.98]{JU10} for the function $k_2$, see \cite{JV2019} for the {\it Back driving L\'evy  process} (BDLP)of $\log \gat$, cf. Section \label{acco} for BDLP's.\\

Recall that the distribution of positive stable r.v.  $\mathbb{S}_\alpha $ is given by \eqref{stap}. As for the r.v. $\gat$,  one can deduce from identity \eqref{logammat},  that the functions, given for $\lambda \geq 0$, by
\begin{eqnarray*}
\phi_{\mathbb{S}_\alpha}(\lambda)&=&-\log \er[e^{-\lambda \mathbb{S}_\alpha}]= t \,\lambda^\alpha=  \int_0^\infty \big(1-e^{-\lambda x} \big)\,\frac{k_3(x)}{x}
\,dx, \quad k_3(x)=\frac{\alpha}{\Gamma (1-\alpha)\,x^{\alpha}}, \\
\Psi_{\mathbb{S}_\alpha}(\lambda)&=& \log \er[(\mathbb{S}_\alpha)^{-\lambda}]=  \int_0^\infty \big(e^{-\lambda x} -1 +\lambda x \big)\, \frac{k_4(x)}{x} \,dx, \quad k_4(x)= \frac{1}{1-e^{-\alpha x}}- \frac{1}{1-e^{-x}}\\
\end{eqnarray*}
are respectively in $\BF$ and $\mathcal{LE}$, and that $k_3, \, k_4$  are non-increasing. By Corollary \ref{char},  it is  immediate  that $\mathbb{S}_\alpha \sim \SD_0(\Ren_+)$ and
$-\log \mathbb{S}_\alpha \sim \ID_-(\mathds{\Ren}) \cap \SD_0(\Ren)$.

\subsection{Main results, new decomposability properties}
After completing this work, we discovered the recent result of  Berg, \c{C}etinkayab \& Karp \cite[Theorem 3.13]{bergkarp} on the function $G_{\alpha ,t}$ given in \eqref{nota}, which is  equivalent to  Theorem \ref{gamXY}. Our Theorem provides a stochastic interpretation of \cite[Theorem 3.13]{bergkarp} and which is shown with a different proof. Furrther, our investigation of the function $G_{\alpha ,t}$  goes beyond by seeking for selfdecomposability properties in Theorem \ref{gself}, Corollary \ref{pidg} and Proposition \ref{maj} below.
\begin{teo} For the function  $G_{\alpha ,t}, \;t>0$, given in \eqref{nota}, we have the equivalences between the assertions  $(i), \;(ii)$ and $(iii)$ in each of the points 1) and 2).
\begin{enumerate}[1)]
\item
    \begin{enumerate}[(i)]
   \item for all $r>0$, the function $(G_{\alpha ,t})^r \in \CM$;

    \item the function $G_{\alpha ,t}$ is  the Laplace transform of some positive infinitely
    divisible  random variable $X_{\alpha ,t}$:
     $$G_{\alpha ,t}(\lambda)  = \er\left[e^{- \lambda \,X_{\alpha ,t}}\right],\quad \lambda\geq 0;$$

    \item  $\alpha \in (0,1)$ and $t \geq 1/2\;$ or $\;\alpha \in (1,\infty)$ and $t < 1/2$.\bigskip
     \end{enumerate}
\item
   \begin{enumerate}[(i)]
   \item the function $(G_{\alpha ,t})^{-r} \in \CM$, for all $r>0$;
    \item the function $(G_{\alpha ,t})^{-1}$ is  the Laplace transform of some positive infinitely divisible random variable $Y_{\alpha ,t}$:
      $$\frac{1}{G_{\alpha ,t}(\lambda)} = \er\left[e^{-\lambda Y_{\alpha ,t}}\right],\quad \lambda\geq 0;$$
    \item $\alpha \in (0,1)$ and $t < 1/2\;$ or $\;\alpha \in (1,\infty)$ and $t \geq 1/2$.
     \end{enumerate}
\end{enumerate}
\label{gamXY}\end{teo}
\begin{rmk}
Equality \eqref{nice}, is equivalent to
$$\er\left[e^{-\lambda  X_{\alpha ,t}}\right]=\er\left[e^{-\alpha \lambda
Y_{1/\alpha,t}}\right]^\alpha , \quad \lambda \geq 0,$$
and this strengthens the fact that $X_{\alpha ,t}$ and $Y_{1/\alpha ,t}$ are concomitantly infinitely divisible.
\end{rmk}
Theorem \ref{gamXY} provides the same result as the one of Lv et al. \cite[Theorem 1.10]{LV} for $G_{\alpha,1}$. Using the fact that the pointwise limit of a sequence of completely monotone functions is also completely monotone, then using Theorem \ref{gamXY} and the fact that the limit
$$G_{\alpha }(\lambda)=\lim_{t\to 0}\, \Gamma(t)^{\alpha-1} \, \alpha^{-\alpha \lambda} \, G_{\alpha ,t}(\lambda),$$
we immediately retrieve the results of  Li \& Chen \cite[Theorem 9]{lichen}
and also of Alzer \& Berg \cite[Theorem 3.5]{Alzer}  on the functions $G_{\alpha}$.  Mehrez  \cite[Theorem 1]{tun} obtained the same result  in the context of the $q$-analog of the function $G_{\alpha}$, i.e. when it is build with the $q$-Gamma function  defined in \eqref{qgam}.
\begin{cor} We have the following results for the function $G_{\alpha}$ given in \eqref{nota}.
\begin{enumerate}[1)]
  \item The function $\lambda \mapsto G_{\alpha}(\lambda)^r$, is completely monotone for all $r>0$, if, and only if, $\;\alpha \in (1,\infty)$;
  \item The function $\lambda \mapsto G_{\alpha}(\lambda)^{-r}$, is completely monotone for all $r>0$, if, and only if, $\;\alpha \in (0,1)$.
\end{enumerate}
\end{cor}
\begin{teo}  Let $t_0=\frac{1}{2}+\frac{1}{2\sqrt{3}}$. It holds that
\begin{enumerate}[1)]
  \item $X_{\alpha ,t} \sim \SD_0(\Ren_+)$ if, and only if, $\alpha \in (0,1)$ and $t \geq t_0$;
  \item $Y_{\alpha ,t} \sim \SD_0(\Ren_+)$ if, and only if, $\alpha \in (1,\infty)$ and $t \geq t_0$.
\end{enumerate}
\label{gself}\end{teo}
In \cite[Lemma 92. Theorem 9.7]{SSV},  it is shown that
$$\ME\subset \BO \subset \ID(\Ren_+)$$
where the Bondesson class $\BO$ is characterized as the smallest class of probability measures on $\oi$ which contains $\ME$, the class of mixture of exponentials distribution, and which is closed under convolutions and vague limits. Actually, the class $\BO$ is one-to-one to the class of {\it complete Bernstein functions}, viz. of functions of the form
\begin{equation}\label{bobo}
\phi(\lambda)= \dd\lambda +\int_{\oi}\frac{\lambda}{\lambda+ x} \frac{\nu(dx)}{x}= \dd\lambda +\int_0^\infty (1-e^{-\lambda x})\, \Lap\nu(x)\,dx ,
\end{equation}
where $\quad  \Lap\nu(x)=\int_{\oi}e^{-x\;u}\nu(du).$ Furthermore, the representation
\begin{equation}\label{bobo1}
\nu(dx)=x \;\eta(x), \quad \mbox{where $\;\eta:\oi \to [0,1],\;$ is measurable},
\end{equation}
is equivalent to say that the Bernstein function $\phi$  is associated to a distribution in $\ME$, cf. \cite[Theorem 9.5]{SSV}. As a consequence of Theorems \ref{gamXY} and \ref{gself}, we obtain the following factorizations for the Gamma distributions.
\begin{cor} Let $\;t>0$, $\underline{\alpha}=(\alpha _1,\alpha _2,\ldots ,\alpha _n) \in(0,1)^n$ such that $\sum_{k=1}^n\alpha _k=1$, $d({\underline{\alpha}}):=\prod_{k=1}^n \alpha_k^{\alpha_k}$ and recall
the random variables $X_{\alpha_k,t}$ and $Y_{\alpha_k,t}$ given by Theorem \ref{gamXY}. Let $\gat$ be a Gamma distributed random variable with shape parameter $t$ and $\mathbb{G}_{1,t},\,\ldots ,\,\mathbb{G}_{n,t}$ denote independent copies of $\gat$. Assuming that the random variables involved in the next factorizations in law are independent, we have the following results.
\begin{enumerate}[1)]
\item  If $t \geq 1/2$ , then we have the  factorization  in law
    \begin{equation} \gat \simdis d({\underline{\alpha}})\;\, \mathbb{G}_{1,t}^{\alpha _1}\,\ldots \,\mathbb{G}_{n,t}^{\alpha _n} \,\;e^{-X_{\underline{\alpha},t}},\quad \mbox{where} \;\; X_{\underline{\alpha},t} := \sum_{k=1}^{n}  X_{\alpha_k,t}.
       \label{idgv}\end{equation}  \label{pidg}

\item   If $\;0<t <1/2$ , then we have the  factorization  in law
    \begin{equation}
      d({\underline{\alpha}})\;\, \mathbb{G}_{1,t}^{\alpha _1}\,\ldots \, \mathbb{G}_{n,t}^{\alpha _n} \simdis \gat  \,\;  e^{- Y_{\underline{\alpha},t}}, \quad \mbox{where}\;\; Y_{\underline{\alpha},t}:=  \sum_{k=1}^n  Y_{\alpha_k,t}.
    \label{idgw}\end{equation}

\item If $\;t\geq t_0=\frac{1}{2}+\frac{1}{\sqrt{12}}$, then  $X_{\underline{\alpha},t}, \; Y_{\underline{\alpha},t}\sim \SD_0(\Ren_+)$.
\item We have $X_{\underline{\alpha},t}, \;Y_{\underline{\alpha},t}\sim \BO,\;$   if, and only if, $\;t\geq 1$.
In this case, their $(1+t)$-fold convolutions, in the sense of \eqref{fold}, have distributions in $\ME.$
\end{enumerate}
\label{pidg}\end{cor}
\begin{rmk} 1) A direct application of Corollary \ref{pidg} is
\begin{equation}
\log \gat \in \ID_-(\mathds{\Ren}) \cap \SD_0(\Ren), \;\, \mbox{if} \;\, t > 1/2 \quad \mbox{and} \quad \log \gat  \in \ID_-(\mathds{\Ren}) \cap \SD_1(\Ren),  \;\, \mbox{if} \,\;t > t_0.
\end{equation}
2) The identities \eqref{idgv} and \eqref{idgw} have an  interpretation in terms of the Gamma function  by taking the Mellin-transforms in both sides: for all $(\alpha _1,\alpha _2,\ldots ,\alpha _n) \in(0,1)^n$,  $\sum_{k=1}^n\alpha _k=1$ and all $r>0$,  we have
  \begin{eqnarray}
  \lambda \mapsto  \er[e^{- \lambda\;X_{\underline{\alpha},t}}]^r &=& \left(\frac{\Gamma(\lambda +t)}{d(\underline{\alpha})^\lambda\prod_{k=1}^n \Gamma(\alpha_k \lambda +t)}\right)^r \in \CM, \quad \mbox{if $t\geq \frac{1}{2}$} \label{rab}\\
  \lambda \mapsto  \er[e^{- \lambda\;X_{\underline{\alpha},t}}]^r &=& \left(\frac{d(\underline{\alpha})^\lambda\prod_{k=1}^n \Gamma(\alpha_k \lambda +t)}{\Gamma(\lambda +t)}\right)^r \in \CM, \quad \mbox{if $t< \frac{1}{2}$}.\label{rab1}
  \end{eqnarray}
The case $t=1$ has been treated by Simon \cite[subsection 3.4]{sim1}. The latter has to be approached to Karp \& Prilepkina's recent result \cite[Theorem 4]{karp}, which states that the function
$$\lambda \mapsto \frac{\prod_{k=1}^{p} \Gamma(A_k \;\lambda + a_k)}{\prod_{j=1}^{q}\Gamma(B_j \lambda + b_j)}, \quad A_k, \; B_j>0, \;a_k,\;b_j\geq 0,$$
is completely monotone if, and only if,
$$\sum_{k=1}^{p} A_k=\sum_{j=1}^{q} B_j, \quad \prod_{k=1}^{p} A_k^{A_k}=\sum_{j=1}^{q} B_j^{B_j}\quad  \mbox{and}  \quad P(u):=\sum_{j=1}^p \frac{e^{-a_j\; u/A_j}}{1-e^{-a_j\; u/A_j}}- \sum_{k=1}^q
\frac{e^{-b_k\; u/B_k}}{1-e^{-b_k\; u/B_k}}\geq 0, \quad \mbox{\it for all}\;  u > 0.$$
As we see, our conditions in Corollary  \ref{pidg}, for $r=1$, are not expressed in terms of the function $P$, since we assume $p=1$, $t=a_k=b_k$ and we require the stronger property of logarithmic complete monotonicity.
\end{rmk}

Using iterates of the shift operators $\Delta_c f(x)= f(x+c) -f(x)$, Akita \& Maejima
\cite[Theorem 1]{makoto} have shown that $\log \gat \sim \SD(\Ren)$ if $t>1/2$.  In their Remark 2, they claimed :
\begin{center}
 ``{\it It is possible to extend to this property to $(t_1, \infty)$ for some $t_1 \in (0,1/2)$}",
\end{center}
that they evaluated, with numerical calculations, by $t_1 \leq 0.152$. The next result and \eqref{t1} in the proof of Proposition \eqref{hats} below provides the value of $t_1$ as the maximum of an explicit elementary function, a quantity that could not be computed  by hand and that we evaluated by Maple.

\begin{prop}  For every $t>0$ and $\alpha \in (0,1)$, we have $\log\gat \sim \ID_-(\mathds{\Ren}) \cap \SD_0(\Ren)$ and the identity in law
\begin{equation}\label{gatg}
\log\gat  \simdis \alpha \log \gat + T_{t,\alpha}.
\end{equation}
Further, $\log\gat \sim \ID_-(\mathds{\Ren}) \cap \SD_1(\Ren)$, i.e. $T_{t,\alpha} \sim \ID_-(\mathds{\Ren}) \cap \SD_0(\Ren)$, if, and only if
$t > t_1\approx 0.151649938034$. The r.v. $T_{t,\alpha}$ corresponds to $\;\log \,J_{\alpha, \alpha t}\;$  in \eqref{privat}.
\label{maj}\end{prop}
\subsection{Comments on the factorizations}
We start by observing the following:
\begin{prop}
If ${(\alpha _k)}_{k\geq 1}$ is a sequence of non increasing positive numbers such that $\sum_{k=1}^\infty\alpha _k=1$ then the sequences $X_{\underline{\alpha},t}$ and $Y_{\underline{\alpha},t}$ of Corollary \ref{pidg} converge  in distribution as $n\to \infty$.
\label{converg}\end{prop}
\medskip

Let $\mathbb{B}_{s,t}, \, s,t>0$, denotes a Beta-distributed random variable with probability density function
$$\frac{1}{B(s,s+t)} \,x^{s-1}(1-x)^{s+t-1}, \quad 0<x<1.$$
It is also worth noticing the following facts.
\begin{enumerate}
  \item  In a private communication, the first author  provided to Bertoin \& Yor their Lemma 1 in \cite{BY}, which shows a factorization close  to \eqref{idgv}. This factorization states that
       if $\mathbb{S}_\alpha,\;0<\alpha <1$, is a standard positive stable r.v. and if $0< \alpha t< s$,  then we have the following factorization in law: with the convention $\mathbb{B}_{t,0}=1$, and
$$\pr(\mathbb{S}_{\alpha,s} \in dx) = \frac{\pr(\mathbb{S}_\alpha \in dx)}{\er[\mathbb{S}_\alpha^{-s}]\; x^{s}},\quad x>0,$$
we have
\begin{equation}\label{privat}
\mathbb{G}_t \simdis \mathbb{G}_s^\alpha \;\mathbb{J}_{\alpha,s}, \qquad \mbox{where}\; \; \mathbb{J}_{\alpha,s}\simdis \frac{ \mathbb{B}_{t,\frac{s}{\alpha}-t} }{  \left(\mathbb{S}_{\alpha,s} \right)^\alpha},\quad \mbox{(on the r.h.s, the r.v.'s are independent)},
\end{equation}
\item Identity  \eqref{idgv}  has to be compared with Gordon's one \cite[Theorem 6]{gordon}: if $p \geq  2$ is an integer,  then
        \begin{equation}\label{g1}
         \frac{\mathbb{G}_{pt}}{p} \simdis \left(\gat  \, \mathbb{G}_{t+\frac{1}{p}} \, \ldots \mathbb{G}_{t+ \frac{p-1}{p}}\right)^{\frac{1}{p}},\quad \mbox{(on the r.h.s, the r.v.'s are assumed to be independent)}.
        \end{equation}
        As an immediate consequence of Gordon's factorization and the one in \eqref{idgv}, we recover a new independent factorization in law for the Beta distributions:  Using the r.v. $X_{\underline{\alpha},t}$ in \eqref{idgv} and the Gamma-Beta algebra, we get
       \begin{equation}\label{g2}
        p\, \mathbb{B}_{t,(p-1)t}\simdis  e^{-X_{\underline{\alpha},t}}   \left(\mathbb{B}_{t,\frac{1}{p}}\,\mathbb{B}_{t,\frac{2}{p}}\ldots \mathbb{B}_{t,\frac{p-1}{p}}\right)^{\frac{1}{p}}, \quad \mbox{(on the r.h.s, the r.v.'s are independent)}.
        \end{equation}
\item Let $0<\alpha \leq 1$ and $ \beta= 1-\alpha$.  Motivated by the selfdecomposability property of the r.v. $\log (\mathbb{S}_\alpha)$,  Pestana, Shanbhag \& Sreehari \cite{PSS}, explored the structure of the r.v. $V_\alpha $ intervening in the so-called Kanter's  identity which involves the exponential and the positive stable distributions:
\begin{equation}\label{pss}
\mathbb{S}_\alpha ^{-\alpha} \; \simdis\; \mathbb{G}_1^\beta  {e^{-V_\alpha }},\quad \mbox{(on the r.h.s, the r.v.'s are assumed to be independent)}.
\end{equation}
Using the well known independent factorisation in law   $\mathbb{G}_1 \simdis  \left( \mathbb{G}_1/\mathbb{S}_\alpha \right)^\alpha$, which is easily justified by
$$\pr (\mathbb{G}_1^{1/\alpha} > \lambda) =\pr (\mathbb{G}_1> \lambda^\alpha) =e^{-\lambda^\alpha} = \er[e^{-\lambda \mathbb{S}_\alpha}]= \pr (\mathbb{G}_1> \lambda \mathbb{S}_\alpha)= \pr \left(\frac{\mathbb{G}_1 }{\mathbb{S}_\alpha}  > \lambda\right),\quad \lambda\geq 0,$$
then taking $n=2$ and $(\alpha_1, \alpha_2)=(\alpha, \beta)$ in \eqref{idgv}, we retrieve:
$$\mathbb{G}_1 \simdis \left( \frac{\mathbb{G}_1}{\mathbb{S}_\alpha} \right)^\alpha \simdis \alpha^{-\alpha}\; \beta^{-\beta}  \;\, \mathbb{G}_{1,1}^{\alpha}\; \,\mathbb{G}_{1,2}^\beta \,\;e^{-X_{\underline{\alpha},1}}\Longrightarrow \mathbb{S}_\alpha ^{-\alpha} \; \simdis\; \alpha^{-\alpha}\; \beta^{-\beta}  \;\, \mathbb{G}_1^\beta  {e^{-V_\alpha }} \Longrightarrow V_\alpha \simdis X_{\underline{\alpha},1} + \alpha \log \alpha +\beta \log \beta,$$
where, on the r.h.s, the r.v.'s are assumed to be independent. Observe that
$$  \er [e^{-\lambda V_{\alpha}}] =\frac{\Gamma (1+\lambda)}{\Gamma (1 +\lambda\alpha )\Gamma (1 +\lambda\beta)},\quad \lambda \geq 0, \quad \mbox{and that}\quad X_{\underline{\alpha},1} \sim \SD_0(\Ren_+) \cap \ME,$$
constitutes an additional information to Kanter's factorization \eqref{pss}.
\end{enumerate}

\medskip
\section{The Multiplicative convolution, the Euler-Mellin differential operator $\Theta$ and $n$-monotone functions}\label{Mellic}
To provide a simple characterization of multiple  selfdecomposability, we give some account on Euler-Mellin differential operator $\Theta=x\,d/dx$ and its relationship with the multiplicative convolution and the concept of $n$-monotone functions.\\

We recall that  the {\it Mellin convolution} (or multiplicative convolution) of two measures $\mu$ and  $\nu$ on $\oi$ is defined by:
$$ \mu \stp \nu (A)=\int_{\oi^2} \II_A(xy)\mu(dx)\nu(dy),\qquad \mbox{if} \; A\;\mbox{is a Borel set of}\; \oi\,.$$
If $\mu$ is absolutely continues with density function $f$, then $\mu \stp \nu$ is the function given by
$$\mu \stp \nu(x)=f \stp \nu(x)= \int_{\oi}f\left(\frac{x}{y}\right) \frac{\nu(dy)}{y}, \quad x>0.$$
\noindent Notice that the integrals above may be infinite when $\mu$ (and/or  $\nu$) is not a finite measure.\\

The {\it Euler-Mellin} differential  $\Theta$  and its discrete  version $\theta_c$, defined by
\begin{equation}\label{operator}
\Theta(g)(x) =xg'(x) \quad \mbox{and}\quad\theta_c(g)(x)=g(x) -g (x/c), \qquad x,\;c>0.
\end{equation}
will be needed in the sequel, for an account on $\Theta$ operators, we suggest  \cite{Butzer}. The iterates of $\Theta$ are denoted by $\Theta^n, \; n\geq 2$ and $\theta_{c_1}\theta_{c_2}\ldots \theta_{c_n}$ denotes the composition of  $\theta_{c_1},\;\theta_{c_2},\ldots, \theta_{c_n}, \;c_1,c_2\ldots ,c_n \in (0,1)$. The following result will be also  needed in the sequel.

\begin{lem} The following is true for every function  $g : (0,\infty) \longrightarrow \Ren$ and  $c \in(0,1)\cup (1,\infty)$:
\begin{enumerate}[1)]
\item  If $g$ is decreasing and $g(0+)<\infty$, then the Frullani integral
$$\int_0^\infty \frac{g(x) -g(x/c)}{x}\, dx  =  \int_0^\infty \frac{\theta_c(g)(x)}{x}\, dx,$$  equals to $ \big(g(\infty)- g(0+)\big)\log c.$\smallskip

\item  Assume further that $g$ is differentiable. Then,
   \begin{enumerate}[(i)]
   \item the function $x\mapsto \theta_c(g)(x) /x $ is increasing (resp. decreasing),  for every fixed $c \in (0,1)$, if, and only if, $g $ is convex (resp. concave);

    \item  the function $x\mapsto  \theta_c(g)(x)$ is increasing (resp. decreasing),  for every fixed $c \in (0,1)$, if, and only if, $x\mapsto  \Theta g(x)$ is decreasing (resp. increasing).
  \end{enumerate}
\end{enumerate}
\label{op}\end{lem}

\begin{rmk}  Since the operators $\Theta$  and $\theta_c$ commute,  it is easily seen that statement $2)(ii)$ in  Lemma \ref{op} extends to $n$-times differentiable functions $g$ via the iterates  of $\Theta$ and the compositions of the operators $\theta_c$. Then, the two following statements are equivalent:
\begin{enumerate}[(i)]
\item  The functions $x\mapsto  \theta_{c_1}\theta_{c_2}\ldots \theta_{c_n}(g)(x)$ are increasing (resp. decreasing), for every $c_1,c_2\ldots ,c_n \in (0,1)$;
\item  The function $x\mapsto  (-1)^n\Theta^n(g)(x)\;$ is increasing (resp. decreasing).
\end{enumerate}
\label{tbg}\end{rmk}
\medskip

It what follows, our aim is to characterise the class $\mathcal{M}_n$  defined in \eqref{Mn}.
\begin{itemize}
\item Williamson \cite{wil} introduced the class of $n$-{\it monotone functions} on $(0,\infty)$ that one can extend  to functions $f:(a,\infty) \to \Ren, \; a\in [-\infty,\infty)$, by reproducing the same arguments of Schilling, Song \& Vondra\v{c}ek, in \cite[Theorem 1.11 and the discussion p.12  given in case  $a=-\infty$]{SSV}. For this, just observe that  $f$ is $n$-monotone on $(a,\infty)$ if, and only, if $f(x+x_0)$ is  $n$-monotone on $(0,\infty)$ for every $x_0>a$. Hence, we will say that $f$ is 1-monotone if $f(x) \geq 0$ for all $x>a$ and if $f$ is non-increasing and right-continuous. The function  $f$ is $n$-monotone on $(a,\infty)$, $n=2,3,\cdots,$  if it is $n-2$ times differentiable,
\begin{equation}\label{fn}
(-1)^j\, f^{(j)}(x)\geq 0, \quad  \mbox{for all} \;x>a,\;\;j= 0, 1,\cdots, n-2,
\end{equation}
and  $(-1)^{n-1}\, f^{(n-2)}$ is non-negative, non-increasing and convex on $(a,\infty)$.

\item Further, with the adaptation of \cite[Theorem 1.11]{SSV}, we can affirm that $f$ is $n$-monotone on $(a,\infty)$ if, and only if, $f$ has the representation
\begin{equation}
f(x)= \cc+\int_{(x,\infty)} (u-x)^{n-1} \, \nu(du) ,\quad x>a\\
\label{tmo}\end{equation}
for some $\cc \geq 0$ and some measure $\nu$ on $(a,\infty)$.
\item Similarly, $f$ is completely monotone on $\Ren$  if and only if, it is $n$-monotone on $(a,\infty)$, for all $n\geq 1$, and all $a<0$; and the latter ensures that $f(x), \; x \in \Ren$, is also represented as in \eqref{cmr}, with some measure $\nu$ on $(0,\infty)$.
\end{itemize}
By \eqref{tmo}, observe that a function $f$  is $n$-monotone on $(0,\infty)$ if, and only, if it is represented by
$$f(x)=\cc+  \big( (1-u)_+^{n-1}  \stp \mu \big)(x) ,\quad x>0.$$
for some $\cc\geq 0$ and some measure $\mu \; \big( \mbox{take $\mu(du) = u^n \; \nu(du)$ in \eqref{tmo} } \big)$.\\

We will now illustrate to which extent the class of $n$-monotone functions  is intimately related to Euler-Mellin's operator. The iterates of the usual  differential operator and of $\Theta$ are linked by these relations: if $g$ is $n$ times differentiable on some interval $I$, then
\begin{equation}
x^n \,g^{(n)}(x) =  \sum_{m=0}^n s(n,m)\,  \Theta^m(g)(x)  \qquad \mbox{and} \qquad\Theta^n(g)(x)= \sum_{m=0}^n S(n,m)\, x^m \,g^{(m)}(x), \quad x\in I,
\label{dert}\end{equation}
where $s(n,m)$ and $S(n,m)$, $0\leq k\leq n,$ denote the Stirling numbers of the first and  second kind, respectively given by the positive numbers
$$s(n,k)= \frac{1}{k!} \frac{d^n}{dx^n}(\log(x+1))^k_{|x=0}  \;\; \mbox{and} \;\; S(n,m)= \frac{1}{m!} \frac{d^n}{dx^n}(e^x-1)^m_{|x=0},$$
cf. \cite{Butzer}.  Notice that $S(n,k)$ is also defined as the number of partitions of the set $\{1, \cdots,n\}$ into exactly $k$ nonempty subsets. It is also known that $s(n,m)=(-1)^{n-m}[_m^n]$, where $[_m^n]$ is the number of permutations in the symmetric group of order $n$  with exactly $k$ cycles. Using  (\ref{dert}), write
$$ x^n\, (-1)^n \,g^{(n)}(x) =   \sum_{m=0}^n [_m^n]\,  (-1)^m\Theta^m(g)(x)$$
and it is immediate that
\begin{equation}\label{dert1}
(-1)^m \Theta^{m}(g)(x)\geq 0, \quad \forall m=0,1,\cdots, n \quad\Longrightarrow \quad  (-1)^n g^{(n)}(x)\geq 0\,.
\end{equation}
Now, assuming that $g(0)=0$ and $g$ is $n$ times differentiable on $\oi$,  consider $h_n=(-1)^n\Theta^n(g)$. If the functions   $h_m:=(-1)^{m-1}\Theta(h_{m-1})= (-1)^{m}\Theta^m(g)$ are such that  $h_m(\infty)=0$ and $h_m(u)/u$ is integrable at $\infty$, for all $m=1,\cdots,n$, then
$$ h_{n-1}(x)=\int_x^\infty \frac{h_n(u)}{u}du= \int_0^\infty h_n\left(\frac{x}{u}\right) \II_{(0,1]}(u) \, \frac{du}{u}= (\II_{(0,1]}\stp h_n)(x), \quad x>0,$$
then, iterating, we get following inversion formulae
$$g=(\II_{(0,1]})^{\stp {n-k}} h_k= (\II_{(0,1]})^{\stp {n}}\stp h_n.$$
Observing that
\begin{equation}\label{power}
\left(\II_{(0,1]}\right)^{\stp n}(x)= \frac{|\log|^n(x)}{n!}\II_{(0,1]}(x),
\end{equation}
one retrieves that $g$ is expressed by
\begin{equation}\label{trans}
g(x)= \int_x^\infty  \log^n\left(\frac{u}{x}\right) \frac{h_n(u)}{u} du\quad x>0.
\end{equation}
Of course, this discussion is informal, but it illustrates the fact that the transform $h_n\mapsto g$ in \eqref{trans} is the inverse  of  the operator $(-1)^n\Theta^n$. Hence,  for "good" functions $g$, for instance if $(-1)^m\Theta^m(g) \geq 0, \; m=1,\ldots ,n,\,$ we will have
$$g(x)= \int_x^\infty  \log^n\left(\frac{u}{x}\right) \, (-1)^n\Theta^n(g)(u) \, \frac{du}{u} \quad x>0.$$
The transform $h\mapsto g$ in \eqref{trans}, is known as the  {\it Hadamard  integral} of order $n+1$, cf. \cite[(18.43) p. 330]{btroj} and \cite{btroj1}.
Further, observing that the difference operators $\Delta_c(f)(x)=f(x+c)-f(x), \;c>0$ are linked to our $\theta_d$ operators by
\begin{equation}\label{thea}
\Delta_c\big(y\mapsto k(e^y)\big)(x)=k(e^c\, e^x)- k(e^x) =   \theta_{e^c}(k)(e^{c+x}), \quad c>0,
\end{equation}
and observing that if $k$ is $m$-times differentiable on $\oi$, then
\begin{equation}\label{dert2}
(-1)^m\frac{d^m}{dx}k(e^x)= (-1)^m \Theta^m(k)(e^x),\quad x\in\Ren,
\end{equation}
and it appears natural to introduce the following class:
\begin{defi} Let $n=1,2,\cdots$. A function $k:\oi\to \oi$, is said to be $\Theta_n$-monotone functions, and we denote $k\in \mathcal{M}_n$, if $x\mapsto k(e^x)$ is $n$-monotone  on $\Ren$. The function $k$ is $\Theta_\infty$-completely monotone if $x\mapsto k(e^x)$ is completely monotone.
\label{melliw}\end{defi}
By \eqref{fn} and \eqref{dert1}, clearly,
\begin{equation}
k\in \mathcal{M}_n \; \Longleftrightarrow \left\{ \begin{array}{lll}
&&k_m:= (-1)^m \Theta^m(k)\geq 0,\quad  \mbox{for all} \;m= 0, 1,\cdots, n-2,\\
&& \mbox{and}\\
&&\mbox{$k_{n-2}$ is non-negative, non-increasing and $x\mapsto k_{n-2}(e^x)$ is convex}
\end{array}\right.
\end{equation}
With these conditions, necessarily $k_n$ is also convex and,  by \eqref{thea}, it becomes clear that
\begin{equation}\label{imp}
k \in \mathcal{M}_n \; \mbox{(resp.  $\mathcal{M}_\infty$)}\;\Longrightarrow \; k \;\mbox{is $n$-monotone on $(0,\infty)$ (resp.  completely monotone).}
\end{equation}
If $k$ is $n$-times differentiable, then by Remark \ref{tbg},
\begin{equation}\label{imp1}
k \in \mathcal{M}_n \; \Longleftrightarrow \; (-1)^n \Theta^n(k).
\end{equation}
The implication \eqref{imp}, left to right, was observed in \cite[Proposition 1.16]{SSV} and the proof there is based on induction on $n$, without formalizing the class $\mathcal{M}_n$.\\

Finally, after our discussion, adapting \cite[Theorem 4.11]{SSV} for $n$-monotone functions, using \eqref{thea} and the $n$-th power multiplicative convolution of the function $\II_{(0,1]}$ given by \eqref{power}, we easily  deduce the full characterization of $\mathcal{M}_n$, i.e.,  the converse to \cite[Proposition 1.16]{SSV} and the analog of \cite[Theorem 4.11]{SSV}  on multiply monotone functions,  for the class $\mathcal{M}_n$:
\begin{lem} Let $k:(0,\infty)\to (0,\infty)$ and $n\geq 1$ . Then, the following assertions are equivalent.
\begin{enumerate}[1)]
\item $k\in \mathcal{M}_n$;
\item $\theta_{c_1}\ldots \theta_{c_m} (k) \geq 0,$ $\quad$ for all $m=1,2,\cdots, n\,$ and $\,c_1,\cdots,c_n \in(0,1)$;
\item $k$ is of the form
\begin{equation}\label{hada}
k(x)=\cc+ \left(\left(\II_{(0,1]}\right)^{\stp \,n-1}\stp \mu\right)(x) = \cc+ \frac{1}{(n-1) !} \int_{(x,\infty)}   \log^{n-1} \left(\frac{y}{x}\right)  \, \frac{\mu(dy)}{y},\quad x>0,
\end{equation}
where $\cc\geq 0$ and the measure $\mu$ is such that  $\int_1^\infty \log^{n-1}(y) \;\frac{\mu(dy)}{y}<\infty$.
\end{enumerate}
Furthermore,  $k\in \mathcal{M}_\infty$ if, and only if, it is represented by
\begin{equation}
k(x)=\cc+    \int_{(0,\infty)} \frac{1}{x^u}   \, \nu(du) ,\quad x>0,
\label{charki} \end{equation}
with some finite measure $\nu$ on $(0,\infty)$.
\label{chark} \end{lem}
\begin{rmk} The integral transform of $\mu$, in \eqref{hada}, is the measure version of the Hadamard  integral of order $n$ in \eqref{trans}.
\end{rmk}

\section{Simple characterization  of multiple selfdecomposable distributions}\label{simple}
By the relations
\begin{eqnarray*}
\Theta(\phi)(\lambda)&=&\frac{d}{dc}\left( e^{-\theta_{1/c}(\phi)(\lambda)}\right) |_{c=1}\\
\theta_{1/c}(\phi)(\lambda)&=&\phi(\lambda)-\phi(c\lambda)=\int_{c\lambda}^\lambda \phi'(t) dt = \int_c^1 \lambda \phi'(s\lambda) ds=  \int_c^1 \Theta(\phi)(s\lambda) \frac{ds}{s}
\end{eqnarray*}
and the fact that $\BF$ is a closed convex cone, a simple proof for the characterization of $\SD_0(\Ren_+)$ could be provided. For instance, see Aguech and Jedidi \cite{aj}, Behme \cite{behme} and  Mai, Schenk \& Scherer \cite{MSS} for the following characterization.
\begin{teo}[\cite{aj, behme, MSS}] Let $X$ be a nonnegative r.v. with cumulant function $\phi(\lambda)=-\log \er[e^{-\lambda X}], \; \lambda \geq 0.\,$ Recall that $\Theta(\phi)$ and $\theta_d(\phi)$ are given by \eqref{operator}.
\begin{enumerate}[1)]
\item If $\theta_c(\phi) \in \BF$ for some $c>1$  or if $\Theta(\phi) \in \BF$, then $\phi \in \BF$.
\item The following assertions are  equivalent.
\begin{enumerate}[(i)]
  \item $X \sim \SD_0(\Ren_+)$;
  \item $\theta_d(\phi)$ is a cumulant function,  for all $\;d>1$;
  \item $\theta_d(\phi) \in \BF$,  for all $\;d>1$;
  \item $\Theta(\phi) \in \BF$;
  \item $\phi  \in \BF$  and is represented by
\begin{equation}\label{brsd}
\phi(\lambda)=\dd \lambda   + \int_0^\infty  \big(1-e^{-\lambda x} \big)\,\frac{k(x)}{x} \,dx, \quad \lambda \geq {\rm 0}
\end{equation}
where $\dd\geq 0$ and $k$ is a  non-increasing function such that $\int_0^1  k(x) \,dx+ \int_1^\infty  \frac{k(x)}{x} \,dx <\infty$.
\end{enumerate}
\end{enumerate}
\label{sd11}\end{teo}
\begin{rmk} Representation  \eqref{brsd} explains identity \eqref{xcy}. Indeed, we have
$$\er[e^{-\lambda Y_c}]= \frac{\er[e^{-\lambda X}]}{\er[e^{-\lambda c X}]}= e^{-\theta_{1/c} (\phi)(\lambda)},$$
and then, after an elementary change of variable, we obtain the representation
\begin{equation}\label{fcl}
\theta_{1/c}(\phi)(\lambda)=  (1-c) \,\dd \,\lambda  + \int_{(0,\infty)}(1-e^{-\lambda x})\frac{\theta_{c}(k)(x)}{x}dx .
\end{equation}
Since $Y\sim \SD_0(\Ren_+)$, then $\theta_{c}(k)(x)/x$ is necessarily   the density function of  a L\'evy measure.
\end{rmk}
Using Proposition \ref{vague} and mimicking the proof of Theorem \ref{sd11}, we can state this proposition without proof:
\begin{prop} We  have $\SD_0(\mathds{\Ren}) \subset  \ID(\mathds{\Ren})$ and the following holds true.
 \begin{enumerate}[1)] \item Let $X\sim \overline{\ID}_-(\mathds{\Ren})$. Then, the following assertions are equivalent.
       \begin{enumerate}[(i)]
      \item $X  \sim  \SD_0(\Ren)$;
      \item The L\'evy-Laplace exponent  $\Psi \in \LE$ associated to $X$ satisfies $\theta_c \Psi \in \LE$ for every $c>0$;
      \item $\Theta(\Psi)$ has the form
       \begin{equation}\label{psisd}
       \Psi(\lambda)=\aaa \lambda +\bb \lambda^2  + \int_0^\infty \big(e^{-\lambda x}-1+\lambda \chi(x)\big)\frac{k(x)}{x}dx ,\quad \lambda\geq 0
        \end{equation}
        for some $\aaa \in \Ren, \,\bb\geq 0$, some truncation function $\chi$ as in \eqref{image} and some non-increasing function $k$.
        \end{enumerate}
\item $X  \sim \ID_-(\mathds{\Ren}) \cap \SD_0(\Ren)$   if, and only if,  its associated L\'evy-Laplace exponent  $\Psi$ has the form  \eqref{psisd} with $\chi(x)=x$, or equivalently $\Theta(\Psi) \in \LE$.
\item Vague closure: $X  \sim  \SD_0(\mathds{\Ren})$ if, and only if, there exist  two sequences of independent r.v's $X_n', \,X_n''\sim \overline{\ID}_-(\mathds{\Ren}) \cap \SD_0(\Ren)$, such that
$X \eqd \lim_{n} (X_n'- X_n'')$.
\end{enumerate}
\label{vaguesd}\end{prop}

\smallskip

To characterize the class $\SD_n(\Ren)$, it is enough to use Proposition \ref{vaguesd} and to focus on the class
$\ID_-(\mathds{\Ren}) \cap \SD_0(\Ren)$ of real-valued selfdecomposable r.v.'s $X$, having a L\'evy-Laplace exponent $\Psi$  of the form
\begin{equation}\label{repsi}
\Psi(\lambda)=\aaa \lambda+ \bb \lambda^2+\int_{(0,\infty)} (e^{-\lambda x} - 1+\lambda x ) \;\frac{k(x)}{x} dx, \quad \lambda \geq 0.
\end{equation}
for some $\aaa \in \Ren$, $\bb\geq 0$ and some non-increasing function $k$  satisfying the  condition \eqref{integg}, i.e.:
$$\int_{(0,1)} x  k(x)  dx + \int_{[1,\infty)} k(x)  dx<\infty.$$
Using the compositions of the Euler-Mellin operators $\Theta,\,\theta_c$ given by \eqref{operator}  and right after, and also the class of $\mathcal{M}_{n+1}$ given by Definition \ref{melliw}, multiple selfdecomposability property is simply characterized as follows:
\begin{cor}[Characterization of multiple selfdecomposability by the Euler-Mellin operator]
Let  $n \geq 1$.
\begin{enumerate}[1)]
\item Assume one of the following
\begin{enumerate}[(a)]
\item $X \sim  \SD_0(\Ren_+)$ and is associated with a Bernstein function $\phi$ and with the  $k$-function given by \eqref{brsd};
\item $X \sim \ID_-(\mathds{\Ren}) \cap \SD_0(\Ren)$ and is associated to a L\'evy-Laplace exponent $\Psi$ to the $k$-function given by \eqref{repsi}.
\end{enumerate}
Then, the following assertions are equivalent.
    \begin{enumerate}[(i)]
    \item $X \sim  \SD_n(\Ren_+)\;   \big(\mbox{\it respectively} \;  X \sim \ID_-(\mathds{\Ren}) \cap \SD_n(\Ren)\big)$;
     \item $k \in \mathcal{M}_{n+1}$, i.e.,  is represented by
      \begin{equation}\label{kk}
      k(x)=\left(\left(\II_{(0,1]}\right)^{\stp \,n-1}\stp \mu\right)(x) =   \frac{1}{(n-1) !} \int_{(x,\infty)}   \log^{n-1} \left(\frac{y}{x}\right)  \, \frac{\mu(dy)}{y},\quad x>0
     \end{equation}
     and $\mu(dx)/x$ is  a measure satisfying \eqref{x11} and  the additional  integrability conditions at infinity
    \begin{equation}\label{loginteg}
    \int_{(1,\infty)} \log^{n+1}(x) \frac{\mu(dx)}{x} <\infty,
    \end{equation}
    (respectively  satisfying only \eqref{x21} if $X\sim \ID_-(\mathds{\Ren}) \cap \SD_0(\Ren)$);
     \item For all $\; m=1,\cdots,n+1$ and $d_1,\cdots,d_{n+1} >1,\;$ the function
     $$(\theta_{d_1}\ldots \theta_{d_m}) (\phi)\quad \big(\mbox{respectively},  \;(\theta_{d_1}\ldots \theta_{d_m}) (\Psi) \big)$$
     is the Bernstein function of some r.v. $Y_m\sim \ID(\Ren_+) \;\;\big($respectively, is the Laplace-exponent  of some r.v. $Y_m\sim \ID_-(\Ren) \big)$;
    \item  For all $m=1, \cdots,n+1$, $\;\Theta^{m}(\phi)   \in \BF\; \; \big($respectively, $\;\Theta^{m}(\Psi)\in \LE \big)$;
    \item  Let $f=\phi$  (respectively, $\,f=\Psi$). There exists $Y\sim \ID(\Ren_+)\;\; \big($respectively, $Y\sim \ID_-(\Ren) \big)$  such that  we have the representation
    \begin{equation}\label{ur}
    f(\lambda)=\frac{1}{(n-1)!}\int_1^\infty g \left(\frac{\lambda}{x}\right)\; \log^{n}(x) \frac{dx}{x},\quad \lambda \geq 0,
    \end{equation}
    where $g$  is the Bernstein function (respectively Laplace exponent) of $Y$.
    \item $\displaystyle X\simdis \int_{\oi} e^{-s^{1/(n+1)}}dZ_s,\,$ where $Z$ is some subordinator (respectively, some spectrally negative L\'evy process) such that
        \begin{equation}\label{jmk}
        \er[\log^{n+1}(1+Z_1)]<\infty \quad \mbox{(respectively $\er[\log^{n+1}(1+|Z_1|)]<\infty$)}.
        \end{equation}
  \end{enumerate}
   Furthermore, in  (v) and (vi), the r.v.'s $Y$ and $Z_1$ have the same distribution.
   \\
\item $X \sim \overline{\ID}_-(\mathds{\Ren}) \cap \SD_n(\Ren)$ if, and only if,  the associated $k$-function and $\mu$-measure satisfy \eqref{kk} and \eqref{loginteg}.
\item We have the same equivalences as in 1) for $n=\infty$, with the following additional conditions on support and integrability on the $\nu$-measure in the representation \eqref{charki}
    of $k\in \mathcal{M}_\infty$:
    \begin{enumerate}[(i)]
    \item  $\nu$ is supported by $(0,1)$, in case $X \in   \SD_\infty(\Ren_+)$ (respectively $(1,2)$, in case $X \in  \ID_-(\Ren)\cap \SD_\infty(\Ren)$), and satisfies
    $$\int_{(0,1)}  \frac{1}{x(1-x)}  \nu(dx)<\infty \; \; (\mbox{respectively}, \;\;
    \int_{(1,2)}  \frac{1}{(x-1)(2-x)} \nu(dx)<\infty).$$
    The latter is equivalent to the representation
    $$\phi(\lambda)=\dd\lambda +\int_{(0,1)} \lambda^x \, \frac{x}{\Gamma(1-x)} \nu(dx)
    \; \; \big(\mbox{respectively}, \;\;
    \Psi(\lambda)= \aaa \lambda +\bb \lambda^2 +\int_{(1,2)} \lambda^x\;\,\frac{x(x-1)}{\Gamma(2-x)} \; \nu(dx)\big),$$
    for some $\dd, \;\bb\geq 0, \; \aaa \in \Ren$.

     \item $\nu$ is supported by $(0,2)$, in case $X \in   \overline{\ID}_-(\Ren)\cap\SD_\infty(\mathds{\Ren})$,
     and satisfies $\displaystyle \int_{(0,2)}  \frac{1}{x(2-x)}  \nu(dx)<\infty.$
    \end{enumerate}
   \end{enumerate}
\label{char} \end{cor}
As a straightforward consequence, we recover the general case:
\begin{cor} Assume  $1\leq n \leq \infty$ and $X \sim  \SD_0(\Ren)$. Then,  $X \sim  \SD_n(\Ren)$ if, and only if $X=X^+-X^-$, where $X^+$ and $ X^-$ are independent and have distributions in $ \overline{\ID}_-(\mathds{\Ren}) \cap \SD_n(\Ren)$.
\label{chacha}\end{cor}
\begin{rmk}
Formula like the above \eqref{kk} appeared in Urbanik  \cite{urb72b, urbanik1}. The stochastic integral representation in 1)(vi) of Corollary \ref{char} was proved in Jurek \cite[Corollary 2.11]{JU2}.
\end{rmk}
\section{Possible extensions of multiple selfdecomposability through  integral stochastic representations}\label{verva}
Observe that the integral stochastic representation \eqref{verv} admits several extensions as noticed by Jurek \& Vervaat \cite{JV1983}. As in the proof of \cite[Proposition 1]{JV1983}, we have the following: for $0<a<b$, a function  of bounded variation $h:(a,b]\to \Ren$, a   monotone function $r:(0,\infty)\to (0,\infty)$, and a L\'evy process $Z(t)$,  the random integral
$$X=I^{h,r}_{(a,b]}(Z):= \int_{(a,b]}h(s)dZ_{r(s)},$$
is well defined and
\begin{equation}
\Phi_{X}(u):= \log \mathbb{E}[e^{\ii \,u \, X}] = \int_{(a,b]}\log \mathbb{E}[e^{\ii\, u \, \,h(s)Z_1}]dr(s),\quad  u \in \Ren.
\end{equation}
By Jurek \cite{JU0},  the Fourier L\'evy-exponents $\Phi_{Z(1)}$ and $ \Phi_{X}$ are linked by
$$ \Phi_{X}(u)=\int_{(a,b]}\Phi_{Z_1}(uh(s))dr(s).$$
The following may be viewed as a particular of the above scheme. For subordinators, Maejima \cite{mama} and  Schilling, Song \& Vondra\v{c}ek \cite{SSV} had the same approach: let $f:(0,A)\to (0,B)$ be a strictly decreasing function and  $Z={(Z_t)}_{t\geq 0}$ be a subordinator with associated Bernstein function $\phi_Z$. Then, the r.v.
\begin{equation}\label{choice}
X = \int_{(0,A)} f(s) dZ_s
\end{equation}
is  a well-defined r.v. on $[0,\infty]$. Reasoning by Riemann approximation of the stochastic integral, \cite[Lemma 10.1]{SSV} provides the Laplace representation
\begin{equation*}
\er[e^{-\lambda X}] = e^{-\phi_X(\lambda)}, \qquad \phi_X(\lambda) = \int_0^A \phi_Z\big(\lambda f(s)\big) ds, \quad \lambda \geq 0.
\end{equation*}
Since $f$ admits a strictly decreasing inverse $F=(0,B)\to (0,A)$, then, by the change of variable $s=G(y):=F(1/y)$, we see that $\phi_X$ takes the form of a Mellin convolution
\begin{equation}\label{fs}
\phi_X(\lambda) = \int_{(1/B,\infty)} \phi_Z\left(\frac{\lambda}{y} \right) dG(y).
\end{equation}
Since $\BF$ is a closed convex cone, then $\phi_X \in \BF$ whenever it is a well defined function on $\Ren_+$. the latter shows that $X\sim \ID(\Ren_+)$. Two particular cases arise for the choice of $f$ in \eqref{choice}:\\

$\bullet$ Taking a continuous nonnegative r.v. $Y$ independent of $Z$ and $f(s)=\pr(Y>s)$, we see that $X$ corresponds to the conditional expectation  $X=\er [Z_Y|Z]$ ,and the integral stochastic representation \eqref{verv} corresponds to the case where $Z$ has the standard exponential distribution. In particular, one can choose
     \begin{equation}\label{choice1}
     f(s)=\er[e^{-sY}]=e^{-\phi_Y(s)}=\pr(\mathbb{G}_1/Y>s), \quad s>0,
     \end{equation}
where $\mathbb{G}_1$ is exponentially distributed and independent of $Z$ and we are also in the previous situation:  $X=\er[Z_{\footnotesize\mathbb{G}_1/Y}\;|\;Z]$.\\

$\bullet$ A  possible extension of the latter  is to take  $f$ of the form
\begin{equation}\label{choice2}
    f=e^{-\phi}\quad \mbox{and $\quad \phi:(0,A) \to (-\log B, \infty)$ differentiable and strictly concave}.
    \end{equation}
Observe that  $\phi$ is necessarily increasing and that the inverse function  $\Psi$, of $\phi$, is a differentiable  increasing and strictly convex function on $\oi$. Making  the change of variable $s=G(y)=\Psi(\log y), \, y>1/B$,  the Bernstein function of $X$ in \eqref{fs} takes the form
\begin{equation}\label{repx}
\phi_X(\lambda) = \int_{1/B}^\infty \phi_Z\left(\frac{\lambda}{y} \right) \Psi'(\log y)\frac{dy}{y}, \quad   \; \lambda \geq 0. \end{equation}
Let $\dd_Z, \,\Pi_Z$ denote  the drift term and the L\'evy measure in the representation \eqref{lev1} of $\phi_Z$.  Since $\Psi'\geq 0$, then Fubini-Tonelli applies in \eqref{repx} and by a change of variable, we obtain the following  Mellin convolution representations:
    \begin{eqnarray}\label{identi}
    \phi_X(\lambda) &=& \dd_Y\, c_{\tiny \Psi} \lambda +\int_0^\infty  (1-e^{-\lambda x}) \frac{l_X(x)}{x}dx,\label{identi} \qquad  c_{\footnotesize \Psi}:=\int_{1/B}^\infty \Psi'(\log y)\frac{dy}{y^2},\\
    \nonumber\\
     l_X(x) &:=&   \int_{(x/B,\infty)} \Psi'\left(\log \frac{u}{x} \right)\Pi_Z(du), \quad x>0. \label{identi11}
    \end{eqnarray}
Assuming \eqref{choice2}, the following is worth to be noticed:
\begin{enumerate}[(a)]
\item $c_{\footnotesize \Psi}=0$ if $B=\infty$.
\item The definiteness of $\phi_X$, i.e., the finiteness of the r.v. $X$ given \eqref{choice}, is guaranteed by the following: if $\dd_Y>0$, then the finiteness of the integral $c_{\footnotesize \Psi}$ is required for the definiteness of $\phi_X$ (i.e. $X$ is a well-defined r.v. on $\Ren_+$).
    Thus, $\phi_X$ is a well-defined function on $\Ren_+$  if, and only if the L\'evy measure $\Pi_Z$ satisfies the additional integrability condition
    \begin{equation}\label{identi1}
    \int_0^\infty (x\wedge 1) \, l_X(x) \,\frac{dx}{x}= \int_{\oi} a_\Psi(u)\; \Pi_Z(du)<\infty, \quad
    \end{equation}
    where
    \begin{equation}\label{aps}
    a_\Psi(u) := \int_{0}^{Bu} (x\wedge 1) \Psi'\left(\log \frac{u}{x} \right) \frac{dx}{x}.
    \end{equation}
    Then,   with some computation, we obtain that $a_\Psi$ is well defined if, and only if  $a_\Psi(1)=c_{\footnotesize \Psi}<\infty$. In this case, we have
     \begin{equation}\label{identi2}
    a_\Psi(u) =   \left\{  \begin{array}{ll}
                   c_{\footnotesize \Psi} u,&\quad {\rm if }\quad 0<Bu \leq 1\\
                  &\\
                    \Psi(\log u)-\Psi(-\log B)+ \displaystyle   u \int_{u}^{\infty}  \Psi'(\log y) \displaystyle \frac{dy}{y^2},&\quad {\rm if }\quad Bu>1.
                  \end{array}     \right. 
    \end{equation}
\item One can release the differentiability assumption on $\phi$ in case $B=1$. Indeed, the assumption  of strict concavity of $\phi$ together with its positivity, insures that $\phi$ is strictly increasing. Thus, almost everywhere, $\phi$ is  differentiable with a strictly decreasing derivative. Hence, almost everywhere, $\Psi$ is  differentiable and $\Psi' > 0$.
\item One can link \eqref{identi2} with  Lemma \ref{lemo} below,  through the observation
\begin{equation}\label{aph}
x_0=\frac{1}{B}  \quad \mbox{and} \quad \rho_h(dy)=\frac{\Psi'(\log y)}{y^2}dy, \quad y>\frac{1}{B} \quad\Longrightarrow\quad \chi_h(u)=a_\Psi(u), \quad u>\frac{1}{B},
\end{equation}
and notice that $Z_s, \;s>0,$ has the L\'evy measure $s\;\Pi_Z$. Then, we immediately obtain the following consequence which constitutes an improvement and a simplification of the statement of Sato's theorem \cite[Theorem 2.6 and Theorem 3.5]{sato2007} in the case of subordinators.
  \begin{cor} Let $X$ be a random variable represented by the stochastic integral in \eqref{choice}.  Then the following assertions are  equivalent.
  \begin{enumerate}[1)]
  \item $X$ is a well defined r.v. on $[0,\infty)$;
  \item $\er[a_\Psi(Z_1)\, \II_{Z_1>1/B}]<\infty$
  \item $\er[a_\Psi(Z_s)\, \II_{Z_s>1/B}]<\infty$ for all $s>0$;
  \item $\int_{(1/B, \infty)} a_\Psi(u)\; \Pi_Z(du)<\infty$.
  \end{enumerate}
  \label{gara} \end{cor}
\item The stochastic integral representation \eqref{choice} for selfdecomposable distributions corresponds to the case where $\varphi(s)=s$, hence $\Psi'\equiv 1$. In all cases, the function $l_X$ is decreasing, because $\Psi'$ is increasing.  Thus,
\begin{equation}\label{hence}
X  \; \mbox{\it is represented by \eqref{choice}, with $f$ as in \eqref{choice2}} \Longrightarrow X\sim\SD_0(\Ren_+).
\end{equation}
\item By Corollary \ref{char}, we know  that if $X\in \SD_n(\Ren_+)$, then $\phi_X$ takes the form
    $$\phi_X(\lambda)=\int_0^\infty  (1-e^{-\lambda x}) \frac{k_X(x)}{x}dx, $$
   such that $k\in \mathcal{M}_{n+1}$ and is represented by \eqref{hada}:
   $$k_X(x) =\frac{1}{n !} \int_{(x,\infty)}\log^{n} \left(\frac{y}{x}\right)\,\frac{\mu_X(dy)}{y},\quad x>0.$$
   Identifying in \eqref{choice1}, \eqref{identi}  and \eqref{identi11}, we necessarily have there
   $$k_X=l_X,\quad\frac{\mu_X(dy)}{y}=\Pi_Z(dy), \quad B=1\quad \mbox{and}\quad\Psi=\Psi_Y \;\,\mbox{s.t.}\;\,\Psi_Y'(x)=\frac{x^n}{n!}, \; x>0.$$
   Since $\Psi_Y(0)=0$, then
   $$\Psi_Y(x)= \frac{x^{n+1}}{(n+1)!}, \; x> 0\quad \mbox{and}\quad \phi_Y(s)=\Psi_Y^{-1}(s)= \big((n+1)! \, s\big)^{\frac{1}{n+1}}, \; s>0.$$
   Observe that $\phi_Y\in \BF$, and most of all, $Y\simdis\mathbb{S}_{1/(n+1)}$  the positive $1/(n+1)$-stable r.v. given by \eqref{stap}. The change of variable, $s\mapsto \Psi_Y(s)=s^{n+1}/(n+1)!$, in  \eqref{choice} with $f$ as in \eqref{choice1},  amounts to change in the clock in the process $Z$ : $Z_s\rightsquigarrow Z_{\Psi_Y(s)}$.  and we retrieve the following result: there exist a subordinator $Z$, such that
     \begin{equation}\label{fraco}
     X \sim \SD_n(\Ren_+) \Longleftrightarrow X=\int_{0}^{\infty}e^{-\big((n+1)!\; s\big)^{1/(n+1)}}\,dZ_s=\int_{0}^{\infty}e^{-s}\,dZ_{\Psi_Y(s)}.
     \end{equation}
          Further, since
     $$ a_\Psi(u):=  \frac{1}{n!}\int_{0}^{u} (x\wedge 1) \left(\log \frac{u}{x} \right)^n  \frac{dx}{x}=\frac{1}{(n+1)!} \big(u\,\II_{0<u \leq 1}+\log^{n+1}(u)\,\II_{u > 1}\big),$$
     then Corollary \ref{gara} insures that $X$ is well defined if, and only if, one of the equivalent conditions holds
     \begin{equation}\label{finiten}
     \int_{[1,\infty)}a_\Psi(u) \Pi_Z(du)<\infty \Longleftrightarrow \er[\log^{n+1}(1+Z_1)]<\infty.
     \end{equation}

    Also observe that the change of variable, $s\mapsto s/(n+1)!$ amounts to change the scale of time in the process $Z$:  $Z_s\rightsquigarrow Y_s:=Z_{s/(n+1)!}$, we simply express $\SD_n$ property by:
    $$X \sim \SD_n(\Ren_+) \Longleftrightarrow X=\int_{0}^{\infty}e^{-s^{1/(n+1)}}\,dY_s,\quad \mbox{for some subordinator $Y$.}$$
    We emphasize that the last discussion could be rephrased identically for $X \sim \overline{\ID}_-(\mathds{\Ren})\cap \SD_n(\Ren)$, and it suggests that a fractional selfdecomposability property of type $\SD_{1/\alpha}(\Ren)$, with $\alpha \in (0,1)$  might be interesting to be formalized, from the stochastic point of view ,since it would be associated to the fractional version of \eqref{fraco}:
    $$X\sim \SD_{1/\alpha}(\Ren) \Longleftrightarrow X=\int_{0}^{\infty}e^{-s^{\alpha}}\,dZ_s=\int_{0}^{\infty}\er[e^{-s \,Y}]\,dZ_s,$$
    where, with our  notations in \eqref{choice1}, $Y\simdis\mathbb{S}_\alpha$  the positive $\alpha$-stable r.v. given by \eqref{stap}. Seeking more involved stochastic and analytical properties for the class of r.v.'s in \eqref{choice} will be, hopefully, the scope of future work.
\item  Theorem 1.1 in \cite{BFJ}  is closely connected to the class of Bernstein functions obtained in form \eqref{repx} and, particularly, investigates whenever the function $x\mapsto x^a \,l_X(x), \; a>-1$ is completely monotone. With our approach, it is immediate that
    $x\mapsto  x^a \, l_X(x)\in \CM$ for arbitrary L\'evy measures $\Pi_Z$, if and only if, in \eqref{identi},
$$B=\infty\quad \mbox{and} \quad \Psi'(u)= \Psi_a'(u) :=e^{a u-e^{-u}},\;  \;u \in \Ren.$$
Assuming the latter, we have the representation
\begin{equation} \label{lx}
x^a \, l_X(x) =   \int_{(0,\infty)} e^{-\frac{x}{u}}\,u^a\, \Pi_Z(du), \quad \Pi_X(dx)=   \,\frac{l_X(x)}{x}\, dx, \quad x>0.
\end{equation}
For $a =1$  (respectively $a =0$), the shape  of the L\'evy measure $\Pi_X$ corresponds to the well known class  $\CB$ complete Bernstein  functions (respectively  $\TB$ of Thorin Bernstein functions), cf. \cite[chapters 6, 8]{SSV}. Choosing $\Psi_a(-\infty)=0$, observe that the inverse function of $f_a$, given by \eqref{choice2}, is provided by
$$\Psi_a(x)=\int_{e^{-x}}^\infty \frac{e^{-u}}{u^{a+1}} du, \quad x\in \Ren,\qquad f_a^{-1}(s)= \Psi_a(-\log s)= \int_{s}^\infty \frac{e^{-u}}{u^{a+1}} du, \quad s>0. $$
Thus,
$$\Psi_a(+\infty)=\int_0^\infty \frac{e^{-u}}{u^{a+1}} du<\infty \Longleftrightarrow a\in (-1,0)\quad \mbox{and} \quad \Psi_a(+\infty)=\Gamma(-a).$$
Further, using the exponential integral  $Ei$ function and \cite[8.214 (1,2)]{grad}, we have
\begin{equation} \label{p0}
\Psi_0(x)=-Ei(-e^{-x}), \quad x>0 \Longrightarrow \lim_{x\to +\infty}x-\Psi_0(x)=\gamma   \Longrightarrow \lim_{x\to \infty} \frac{\Psi_0(x)}{x}=1,
\end{equation}
where $\gamma$ is the Euler-Mascheroni constant. If $a>0$
\begin{equation} \label{p1}
\Psi_a(x)=\int_{e^{-x}}^1 \frac{e^{-u}}{u^{a+1}} du+\Psi_a(0)\geq  \frac{e^{ax}-1}{a\, e}, \quad x>0\Longrightarrow \lim_{x\to +\infty} \frac{\Psi_a(x)}{x}=+\infty.
\end{equation}
\item For $a\geq 0$, additional investigations on the inverse function $\phi_a$ of $\Psi_a$ are feasible.
   Like $f_a$, the function $\phi_a$ is also not explicit, but it certainly increases and satisfies
   $$\phi_a(0+)= -\infty, \quad \phi_a(\infty)=+\infty \quad \mbox{and}\quad\phi_a(s_a)= 0,\quad \mbox{where} \;\,s_a :=\Psi_a(0)=\int_1^\infty \frac{e^{-u}}{u^{a+1}} du, \quad a\geq 0.$$
   Further, since
   $$\phi_a'(x)= \frac{1}{\Psi_a'(\phi_a(x))}=e^{-a\phi_a(x)+e^{-\phi_a(x)}}=h_a(\phi_a(x)), \quad x>0,\qquad h_a(u):= e^{-au+e^{-u}}=\sum_{n=0}^{\infty}e^{-(a+n)u},$$
   Thus, the function $x\mapsto \hat{\phi} _a(\lambda):={\phi_a}(s_a+\lambda), \; \lambda\geq 0,$  satisfies
   $$\hat{\phi} _a(0)=0, \quad\hat{\phi}_a'(0)=\frac{1}{\Psi_a'(0)}=e\quad  \mbox{and} \quad (\hat{\phi} _a)^{-1}(x)=: \hat{\Psi}_a(x)=\Psi_a(x) -\Psi_a(0)=\int_{e^{-x}}^1 \frac{e^{-u}}{u^{a+1}} du, \;\; x\geq 0,$$
   Observing that $h_a(u), \; u>0$, is a completely monotone function, then, using Faa di Bruno's formula and induction, we retrieve that ${\phi_a'}_{|(s_a, \infty)}\in \CM$, hence we have
   $$ \hat{\phi}_a'= e^{-a\,\hat{\phi}_a+e^{-\hat{\phi}_a}}  \quad \mbox{and}\quad \hat{\phi} _a\in \BF.$$
   Further, by \eqref{p0} and \eqref{p1}, we have
   $$\dd_a:=\lim_{\lambda \to \infty} \frac{\hat{\phi}_a(\lambda)}{\lambda}= 1, \;\mbox{if} \; a=0 \;\;\mbox{and} \; 0\;\mbox{otherwise}.$$
   We deduce that $\hat{\phi}_a$ has the drift term $\dd_a$ and  is associated to  a L\'evy measure $\Pi_a$ and to a subordinator $Y_a={(Y_{a,t})}_{t\geq 0}$ through the representations:
   $$\hat{\phi} _a(\lambda) =\dd_a \lambda+\int_{\oi}(1-e^{-\lambda x}) \Pi_a(dx), \quad \lambda\geq 0,$$
   and
   \begin{equation} \label{curu}
   \hat{\phi}' _a(\lambda)=\dd_a  +\int_{\oi} e^{-\lambda x}\, x\, \Pi_a(dx) = e^{-a\hat{\phi} _a(\lambda)+e^{-\hat{\phi} _a(\lambda)}}= \sum_{n=0}^{\infty}\frac{e^{-(a+n)\hat{\phi} _a(\lambda)}}{n!}= \sum_{n=0}^{\infty}\frac{\er[e^{- \lambda Y_{a, a+n}}]}{n!}.
   \end{equation}
   Thus, letting $\lambda$ to 0,  we get
   $$\lim_{\lambda \to +\infty}\hat{\phi}' _a(\lambda)=\dd_0=1 = \sum_{n=0}^{\infty}\frac{1}{n!}\pr( Y_{0,n}=0)=1+\sum_{n=1}^{\infty}\frac{1}{n!}\pr(Y_{0,n}=0),$$
   which entails
   $$\sum_{n=1}^{\infty}\frac{1}{n!}\pr(Y_{0,n}=0)=0\Longrightarrow\pr(Y_{a,n}=0)=0, \quad \forall n=0, 1, \ldots \Longrightarrow\pr(Y_{a,t}=0)=0, \quad \forall t\geq 0.$$
   Then, applying the Laplace inversion in \eqref{curu}, we get
   $$ x \Pi_0(dx)=\sum_{n=1}^{\infty}\frac{1}{n!}\pr(Y_{0,n}\in dx),\quad x>0.$$
   Similarly,  if $a>0$,
   \begin{eqnarray*}
   \dd_a = 0&=&\sum_{n=0}^{\infty}\frac{1}{n!}\pr(Y_{a,a+n}=0)\Longrightarrow \quad\pr(Y_{a,t}=0)=0, \; \forall t\geq 0,\\
   x \Pi_a(dx)&=&\sum_{n=0}^{\infty}\frac{1}{n!}\pr(Y_{a, a+n}\in dx),\quad x>0.
   \end{eqnarray*}
   If  $N$  is a r.v. with Poisson distribution, with rate parameter equal to 1, independent of the subordinator $Y_a$, then, by \eqref{curu}, we have the closed forms obtained with the help of the subordinated r.v. $Y_{a, a+N}$: for all $\lambda \geq 0$,
   \begin{equation}\label{rem1}
   \hat{\phi} _0'(\lambda) = 1+ e\sum_{n=1}^{\infty}\frac{e^{-1}}{n!}\er[e^{- \lambda Y_{0,n}}]=
   1+ e \er[e^{- \lambda Y_{0,N}}\II_{N\geq 1}] \Longrightarrow  \hat{\phi} _a(\lambda) = \lambda +e\;\er\left[\frac{1-e^{- \lambda Y_{0,N}}}{Y_{0,N}}\II_{N\geq 1}\right],
   \end{equation}
   and if $a>0$,
   \begin{equation}\label{rem2}
   \hat{\phi} _a'(\lambda) = e\,\sum_{n=0}^{\infty}\frac{e^{-1}}{n!}\er[e^{- \lambda Y_{a, a+n}}]=
   e\, \er[e^{- \lambda Y_{a, a+N}}]\Longrightarrow \hat{\phi} _a(\lambda) = e\;\er\left[\frac{1-e^{- \lambda Y_{a, a+N}}}{Z_{a, a+N}}\right].
   \end{equation}
   Additionally, by \eqref{curu} again, we have that for all $a\geq 0, \,\lambda \geq 0$ and $t>0$,
  \begin{equation}
   \er[Y_{a,t}e^{-\lambda Y_{a,t}}]=t\hat{\phi}_0'(\lambda)e^{-t\hat{\phi}_a(\lambda)}=t\sum_{n=0}^{\infty} \frac{\er[e^{- \lambda Y_{a,a+t+n}}]}{n!}= e\,t\, \er[e^{- \lambda Y_{a,a+t+N}}].
   \label{curuu}\end{equation}
   If $\overline{Y}_{a,t}$ is a version of the size biased distribution of $Y_{a,t}$, defined for $t>0$, by
    \begin{equation}
   \overline{Y}_{a,t}  \;\simdis \;  \frac{x\,\pr(Y_{a,t}\in dx)}{\er[Y_{a,t}]}= \frac{x\,\pr(Y_{a,t}\in dx)}{e\,t},
   \label{curu1}\end{equation}
   then, from \eqref{curuu} we get for $a\geq 0$:
   $$\er[e^{-\lambda \overline{Y}_{a,t}}] = \er [e^{-\lambda   Y_{a,a+t+N}}] , \quad \forall \lambda \geq 0 \Longleftrightarrow \overline{Y}_{a,t} \;\simdis\; Y_{a,a+t+N}\simdis\; Y_{a,t} +Y_{a,a+N} , \quad \forall t> 0,$$
   where in the last identity, the r.v.'s $Y_{a,t} $ and $Y_{a,a+N}$ are assumed to be independent.\\

   Finally, by \eqref{choice}, we conclude that  a  nonnegative infinitely divisible r.v. $X$  has a Bernstein function $\phi$  with L\'evy measure of the form  \eqref{lx} and $a\geq 0$, if, and only if, for  some subordinator $Z$, it has the following stochastic  representation
   $$X=\int_{(0,\infty)} f(s) dZ_s,\qquad f(s)=e^{- \phi_a(s)},\quad s\geq 0$$
   and recall that $f$ coincides with
   $$f(s)= e^{-\hat{\phi} _a(s-s_a)}=\er[e^{(s_a-s) Y_{a,1}}],\quad \mbox{if}\quad s\geq s_a.$$
   In \cite[Chapter 10]{SSV}, several forms of the function $f$ were investigated in the purpose of characterizing the induced classes of distribution. The cases $a=0$ and $a=1$ leading to the classes $\TB$ and $\CB$ were also investigated there but no closed-form was proposed for the corresponding function $f$. For this reason, and because of the remarkable relations \eqref{rem1}, \eqref{rem2}, and \eqref{curu1}, the distribution of $Y_{a,1}$ would benefit from being investigated  in more detail.
\end{enumerate}
\section{The proofs}\label{proo}
\begin{proofof}{\bf Proposition \ref{vague}.}  1) Consider the L\'evy-Laplace exponent $\Psi$ of $X\sim \overline{\ID}_-(\mathds{\Ren})$. The  L\'evy measure  $\pi$  gives no mass to $(0, \infty)$ and
its image $\Pi$  given by \eqref{image} satisfies \eqref{x21}.  Observe that the truncated L\'evy measures $\Pi_n:= \Pi_{|(0,n]}, \, n\in \NN$, satisfies \eqref{integg}, so that the quantities
$$\aaa_n:=  \aaa + \int_{(0,n]} \big(\chi(x)-x\big) \, \Pi(dx)\quad \mbox{are finite},$$
and each of the  functions,
\begin{equation}\label{lolo}
\Psi_n(\lambda):=  \aaa_n \lambda \,  + \bb\, \lambda^2 + \int_{(0,\infty)}(e^{-\lambda x}-1 +\lambda \, x) \;\Pi_n(dx),
\end{equation}
belongs to $\LE$, hence, is associated with a r.v. $X_n\sim \ID_-(\Ren)$. Since $\Psi=\lim_{n}\Psi_n$, deduce that $X \eqd \lim_{n} X_n$.
The converse is obvious because $\ID_-(\Ren) \subset \overline{\ID}_-(\mathds{\Ren})$. Statement 2) is an immediate consequence of 1).
\end{proofof}
\subsection{Additional results, preliminaries for the proofs of the main results}
The following technical result will be used in the proof of Lemma \ref{hatam}.
\begin{prop} Let $t\geq 0,\;\alpha \in (0,1)\cup (1,\infty)$ and define for $u >0,$
\begin{eqnarray}
e_{\alpha ,t}(u)&:=& \frac{e^{-t u}}{1-e^{-u}}\label{eat}\\
g_{\alpha ,t}(u)&:=& \alpha\, e_{\alpha ,t}(u) - e_{\alpha ,t}(u/\alpha),\label{gat}\\
h_{\alpha ,t}(u)&:=&  e_{\alpha ,t}(u) - e_{\alpha ,t}(u/\alpha).\label{hat}
\end{eqnarray}
\begin{enumerate}[1)]
\item The following holds for $g_{\alpha ,t}$:
\begin{enumerate}[(i)]
\item  The integral $\int_0^\infty g_{\alpha ,t}(u)du$ is finite if, and only if, $t>0$. In this case, the integral equals to $-\alpha \log \alpha $.

\item  The function $g_{\alpha ,t}(u)$ is nonnegative  (respectively nonpositive) in $u$, for every fixed
$\alpha <1$ (respectively $\alpha >1$), if, and only if, $t \geq \frac{1}{2}$.

\item  The function $g_{\alpha ,t}(u)$ is decreasing (respectively increasing) in $u$, for every fixed $\alpha \in(0,1)$ (respectively $\alpha >1$), if, and only if, $t\geq t_0 =\frac{1}{2}+\frac{1}{2\sqrt{3}}$.
\end{enumerate}
\item The function  $h_{\alpha ,t}(u)$  is  decreasing in $u>0$, for every fixed $\alpha \in(0,1)$ if, and only if, $t\geq t_1 \simeq0.151463487259.$
\end{enumerate}
\label{hats}\end{prop}
\begin{proofof}{\bf Proposition \ref{hats}.} 1)   Since $g_{\alpha ,t}(u)= -\alpha g_{1/\alpha,t}(u/\alpha )$, it is enough to consider the case $\alpha <1$. To understand  the rest of the proof of 1), notice the following expression:
\begin{equation}
g_{\alpha ,t}(\alpha u)= \frac{1}{u}\left[ \frac{\alpha u e^{-t \alpha u}}{1-e^{- \alpha u}}- \frac{u e^{-tu } }{ 1-e^{-u}}\right]= - \frac{\theta_{1/\alpha}(g_t)(u)}{u}, \quad g_t(u)=\frac{u e^{-tu } }{ 1-e^{-u}}.
\end{equation}
$1)(i)$ The necessary and sufficient condition $t>0$ is trivial. The value of the integral is a direct application of 1) in Lemma \ref{op}.

\noindent $1)(ii)$  The necessary part stems from $\lim_{0+}g_{\alpha ,t}(u)=(1-\alpha )(t-\frac{1}{2})$. For the sufficient part, just notice that $\lim_{+\infty}g_{\alpha ,t}(u)=0$ and  write
$$g_t(u)= \frac{u e^{-tu }}{1-e^{-u}} = \frac{u e^{-u/2 }}{1-e^{-u}} e^{-u(t-1/2) }= \frac{u/2}{\sinh(u/2)} e^{-u(t-1/2) }, $$
in order to obtain that $g_t$ is decreasing and to conclude that $g_{\alpha ,t}$ is nonnegative.\\
$1)(iii)$ According to Lemma \ref{op} 2)$(i)$, we only need to check whenever $g_t$ is convex. Standard calculations lead to
$$g_t''(x)= \frac{g_t (x)}{x(e^x-1)^2} P(x,t), \quad x,\,t>0,$$
where
$$P(x,t)= x(e^x-1)^2 t^2 -2t(e^x-1)(e^x-1-x) + x(e^x+1)-2(e^x-1),$$
and it is clear that the convexity of $g_t$ is equivalent to the positivity of the function $P$. Since the discriminant $\Delta_P(x)$ of the polynomial $t\mapsto P(x,.t)$  equals to $$\Delta_P(x)= 4 (e^x-1)^2\left[(e^x-1)^2-x^2e^x \right]=16  e^x  (e^x-1)^2 \left[\sinh(\frac{x}{2})^2 -(\frac{x}{2})^2 \right],$$
and is positive, we recover two positive roots:
$$t_{\pm}(x)= \frac{1}{x(e^x-1)}\left[(e^x-1-x) \pm \sqrt{(e^x-1)^2 -x^2 e^x}\right]= \left[\frac{1}{x} -\frac{1}{e^x-1}  \pm \sqrt{\frac{1}{x^2} -\frac{e^x}{(e^x-1)^2} }\right].$$
Since
$$\left(\frac{1}{x} -\frac{1}{e^x-1}\right)'=- \frac{(e^x-1)^2-x^2 e^x}{x^2(e^x-1)^2},\quad
\left(\frac{1}{x^2} -\frac{e^x}{(e^x-1)^2}\right)' = \frac{e^x (e^x+1)}{(e^x-1)^3}-\frac{2}{x^3},\quad x>0,$$
we  see that  the functions
$$x\mapsto \frac{1}{x} -\frac{1}{e^x-1}\quad\mbox{and}\quad x\mapsto \frac{1}{x^2}-\frac{e^x}{(e^x-1)^2},$$
are both decreasing, and then so is  $x\mapsto t_{+}(x)$. By expansion near $0$, we obtain
$$t_0:=\max_{x>0} t_{+}(x)= t_{+}(0+)=\frac{1}{2}+\frac{1}{\sqrt{12}}.$$
Since $\min_{x>0} \,t_{-}(x)=0$, we deduce that $P(x,t)$ is positive for every $x>0$ if, and only if, $ t\geq t_0$.\\
2) Observe that
\begin{equation} \label{hta}
h_{\alpha ,t}(u)= \theta_{\alpha}(e_t)(u), \quad  e_t(u)=e^{-t u}\,e_0(u), \quad u>0,
\end{equation}
and recall that we are looking for the range of $t$ for which  we have $h_{\alpha ,t}$ is decreasing. By Lemma \ref{chark}, we have
\begin{equation}\label{etta}
\theta_{c_1}\theta_{c_2}(e_t)\geq 0, \quad \forall \; c_1,\,c_2 \in (0,1)\Longleftrightarrow   e_t\in \mathcal{M}_2 \Longleftrightarrow (-1)^n \Theta^n(e_t) \geq 0, \quad \;\mbox{for} \; n=1,2.
\end{equation}
Using the fact that $e_t$ is decreasing, and using \eqref{dert}, it only remains to check whenever
$$\Theta^2(e_t)(x)=x^2\,e_t''(x)+x\,e_t'(x)=
x\,e^{-t u}\big[x\,e_0(x)\,t^2 -\big(e_0(x)+2x e_0'(x)\big)\,t+e_0'(x)+x\,e_0''(x)\big]\geq 0.$$
Since $e_0'=e_0(1-e_0)$ and $e_0''=e_0(1-e_0)(1-2e_0)$, we find that
$$\Theta^2(e_t)(x)\geq 0 \Longleftrightarrow Q(x,t):= xt^2 -\big(1+2x(1-e_0(x)\big)\,t+
\big(1-e_0(x)\big)\big(1+x (1-2e_0(x)\big)\geq 0 .$$
As in 1),  the discriminant $\Delta_Q(x)$ of the polynomial $t\mapsto Q(x,t)$, is given by
$$\Delta_Q(x)=1-4x^2 e_0(x)  \big(e_0(x)-1\big)=1-\frac{4x^2 e^x}{(e^x-1)^2}=\left(1-\frac{x}{\sinh(x/2)}\right)\left(1+\frac{x}{\sinh(x/2)}\right)$$
Paradoxically, things are not as smooth as for $g_{\alpha ,t}$: $\Delta_Q(x)$
is positive if, and only if, $x$ is bigger than some value $x_0$ which can not be expressed by hand, but could be evaluated by Maple with the value $x_0\sim 4.35463796993$. Hence, for $x>x_0$,
we recover two positive roots
$$t_{\pm}(x)=\frac{1}{2x}\left(1-\frac{2x}{e^x-1}\pm \sqrt{1-\frac{4x^2 e^x}{(e^x-1)^2}}\right)$$
and  finally, using Maple again, we get
\begin{equation}\label{t1}
Q(x,t) \geq 0, \quad \forall \; x>0 \Longleftrightarrow t\geq t_1=\max_{x>x_0}\,t_+(x)\simeq 0.151463487259.
\end{equation}
\end{proofof}
\medskip

Let $t>0$,  ${\alpha}=(\alpha _1,\alpha _2,\ldots ,\alpha _n)\in (0,1)^n$. Recall $e_t$ is given by \eqref{eat} and let the functions
\begin{equation} \label{gta}
h_{\underline{\alpha},t}(u)= e_t(u)- \sum_{k=1}^n e_t(u/\alpha _k) ,\quad u >0.
\end{equation}

Alzer \& Berg \cite[Lemma 2.7]{Alzer} provided a Petrovi\'{c}-type inequality for the function $e_1$:
$$\frac{n-1}{2}+\sum_{k=1}^n e_1(u/\alpha_k)- e_1(u/\sum_{k=1}^n\alpha_k)\geq 0, \quad   \mbox{for all $u>0$}.$$

We shall provide additional information for the function $h_{\underline{\alpha},t}$ in Lemma \ref{hatam} below. For this purpose, we need some preliminary results. We denote by $\delta_a$  the Dirac measure in $a$ and by $\lfloor\,\rfloor$ and $\{\,\}$ the  integer and the fractional part functions, respectively.
\begin{lem} Let  $\;t>0$, $\underline{\alpha}=(\alpha _1,\alpha _2,\ldots ,\alpha _n)\in (0,1)^n$ such that $\sum_{k=1}^n\alpha _k=1$ and let $\mu_{\underline{\alpha} ,t}$ and $m_{\underline{\alpha} ,t}$ be the signed measure and the function defined by
$$\mu_{\underline{\alpha} ,t}(dx):=\sum_{i=0}^\infty \big(\delta_{i+t}(dx) - \sum_{k=1}^n\delta_{\frac{i+t}{\alpha_k}}(dx)\big) \quad \mbox{and}\quad m_{\underline{\alpha} ,t}(x):= \mu_{\underline{\alpha} ,t}([0,x)),\; x\geq t.$$
Then,  $m_{\underline{\alpha} ,t}$  is positive for every fixed $\underline{\alpha}$ if, and only if, $t\geq 1$. In this case, we have the control
\begin{equation}
0\leq m_{\underline{\alpha} ,t}(x) \leq 1+t.
\label{rt}\end{equation}
\label{debloc}\end{lem}
\begin{proofp} We treat only the case $n=2$, the case $n\geq 3$ is treated similarly.  By symmetry  of the problem, we may assume that $0<\alpha<\beta=1-\alpha<1$ and take $\underline{\alpha}=(\alpha, \beta)$. Since for every $\gamma >0$,
$$\sum_{i=0}^\infty\delta_{\frac{i+t}{\gamma}}\big([0,x)\big)= \lfloor \gamma x+ 1-t \rfloor\,\II_{( x\geq t/\gamma)},$$
then
\begin{eqnarray}
m_{\underline{\alpha} ,t}(x) =\mu_{\underline{\alpha} ,t}([0,x)) &=& \lfloor x+ 1-t\rfloor\,\II_{(x\geq t)} -\lfloor \alpha x+ 1-t\rfloor\,\II_{( x \geq t/\alpha)}-\lfloor \beta x+ 1-t\rfloor\,\II_{(x\geq t/\beta)} \nonumber \\
&=&\left\{
\begin{array}{ll}
  0& {\rm if }\;\; 0<x < t\\
\lfloor x-t\rfloor +1 & {\rm if }\;\; t\leq x < t/\beta\\
\lfloor x-t\rfloor-\lfloor\beta x-t\rfloor& {\rm if }\;\; t/\beta \leq x < t/\alpha\\
\lfloor x -t\rfloor-\lfloor\alpha x-t\rfloor-\lfloor\beta x-t\rfloor-1& {\rm if }\;\; t/\alpha \leq x.
\end{array} \right. \label{carc}
\end{eqnarray}
\noindent  1) Assume $t\geq 1$.   Using the fact that $y -1 < \lfloor y\rfloor \leq y$,  for all $y \in \Ren,$  the upper bound in \eqref{rt} easily follows.  Because the integer part function  is increasing, we just need to study the positivity of the function $m_{\underline{\alpha} ,t}$  restricted to $(t/\alpha, \infty)$.  Since $ \lfloor y+z\rfloor  \geq  \lfloor y\rfloor + \lfloor z\rfloor,\;\forall y, z \in \Ren$, we deduce that, for every  $x > t/\alpha$,
$$m_{\underline{\alpha} ,t}(x)= \lfloor x+1-t\rfloor-\lfloor\alpha x+  1-t + \beta x+  1-t \rfloor \\
= \lfloor x+1-t\rfloor-\lfloor x+ 2(1-t)\rfloor \geq 0.$$
2) If  $t< 1$, then the function $m_{\underline{\alpha} ,t}(x)$ could take negative values if $x > t/\alpha$. Indeed, choose  $x$ such that for some positive integers $k$ and $l$, we have
$$k\leq \alpha x -t<k+\frac{1-t}{2}\quad \mbox{and}\quad l\leq \alpha x -t<l+\frac{1-t}{2},$$
so that
$$\{ \alpha x -t\}<\frac{1-t}{2}\quad \mbox{and}\quad \{\alpha x -t\}<\frac{1-t}{2}.$$
and then
\begin{eqnarray*}
m_{\underline{\alpha} ,t}(x)&=& \lfloor x-t\rfloor-\lfloor\alpha x-t\rfloor - \lfloor\beta x -t\rfloor -1=t+ \{\alpha x-t\} + \{\beta x -t\} -1-\{ x-t\} \\
&\leq &t+ \{\alpha x-t\} + \{\beta x -t\} -1 <0.
\end{eqnarray*}
\end{proofp}

This Lemma will be used in the proof of Corollary \ref{pidg}.
\begin{lem} Let $t_0$ be the universal constant of Theorem \ref{gself}. The function $h_{\underline{\alpha},t}$  defined by  \eqref{gta} with $\sum_{k=1}^n\alpha _k=1$ satisfies the following:
\begin{enumerate}[1)]
\item If $t\geq 1/2$, then  $ h_{\underline{\alpha},t}(u)>0,\;$  for all $u >0$.
\item If $t\geq t_0$, then   $ h_{\underline{\alpha},t}$\; is decreasing.
\item The function $u\mapsto h_{\underline{\alpha},t}(u)/u\;$ is completely monotone if, and only if, $t\geq 1$. In this case, we have the representation
$$h_{\underline{\alpha},t}(u)  = (1+t) \int_0^\infty e^{-ux } \,\eta(x)\, dx, \quad u>0,$$
for some measurable function $\eta : (0,\infty)\to [0,1]$.
\end{enumerate}
\label{hatam}\end{lem}
\begin{proofp} 1) and 2) Just notice the expression  $h_{\underline{\alpha},t}(u)= \sum_{k=1}^n g_{\alpha _k,t}$, with $g_{\alpha _k,t}$ given by \eqref{gat},  and apply Proposition \ref{hats}.\\
3)  Expanding  the terms in  $\ga$, and using the signed  measure and the positive function $\mu_{\underline{\alpha} ,t}$ and $m_{\underline{\alpha},t}$ of Lemma \ref{debloc}, obtain the expression
$$ \frac{h_{\underline{\alpha},t} (u)}{u}= \frac{1}{u}\sum_{i=0}^\infty \big(e^{-(i+t)u}- \sum_{k=1}^n e^{-\frac{(i+t)}{\alpha_k }u}\big)=\frac{1}{u}\int_{[t,\infty)} e^{-ux}\mu_{\underline{\alpha} ,t}(dx) \\
=\int_t^\infty  e^{-ux} m_{\underline{\alpha} ,t}(x)\, dx, \quad u>0,$$
and conclude with the nonnegativity of $m_{\underline{\alpha} ,t}$.
\end{proofp}
\subsection{Linking the integrability of infinitely divisible distributions with their  L\'evy measure} \label{lik}
Let $\mathcal{C}$ be the class of functions $h:[0,\infty) \to [0,\infty)$, differentiable, such that
\begin{equation}\label{coca}
 \lim_{x\to \infty}h'(x)=0\quad \mbox{there exists $ x_0 \geq 0$ s.t. $h$ is concave on $[x_0,\infty)$}.
\end{equation}
Note that for such functions, there exists a finite positive measure $\rho_h$ on $[x_0,\infty)$ such that  the derivative of $h$ is  represented by
\begin{equation}\label{coca1}
h'(x)=\rho_h \big([x,\infty)\big), \quad  x\geq x_0.
\end{equation}

The following lemma  gives an interpretation of the integrability condition \eqref{identi1} and constitutes a variant of \cite[Theorem 2]{JUS}, which is stated with sub-multiplicative functions $h$.
\begin{lem}
Let $Z$ be nonnegative r.v. with cumulant function $\phi(\lambda)=-\log \er[e^{-\lambda Z}]$ and  let a function $h$ in the class $\mathcal{C}$ defined by \eqref{coca} and associated to the pair $(x_0,\rho_h)$ by \eqref{coca1}.
     \begin{enumerate}[1)]
     \item     We have the equivalence:
     \begin{equation}\label{ffin3}
    \er[h(Z)]<\infty \Longleftrightarrow \int_{x_0}^{\infty}  \min\big(1,\phi(1/x)\big)\;x\; \rho_h(dx)  <\infty;
     \end{equation}
    \item Assume $Z\sim \ID(\Ren_+)$ and has characteristics $(\dd, \,\Pi)$ in the representation \eqref{lev1} of its Bernstein function $\phi$, then the following assertions are equivalent.
        \medskip

        \begin{enumerate}[(i)]
           \item $\er[h(Z)]<\infty$;
           \smallskip
           \item $\int_{[x_0,\infty)}  \chi_h(u)\, \Pi(du) <\infty$, where
             \begin{equation} \label{chi}
             \chi_h(u)= \int_{[x_0,\infty)}(x\wedge u)\;\rho_h(dx)= \left\{ \begin{array}{lll}
             && x_0\;h'(x_0) \;u, \quad \mbox{if}\;\; 0<u<x_0,\\
             && \\
             &&h(u)-h(x_0)+ x_0\;h'(x_0),\quad  \mbox{if}\;\;u \geq x_0.
             \end{array}\right.
             \end{equation}
          \item $\int_{[x_0,\infty)}  \big(h(u)-h(x_0)\big)\; \Pi(du) <\infty$.
        \end{enumerate}
   \end{enumerate}
  \label{lemo}   \end{lem}
\begin{proofof}  1)  Since $h'(\infty)=0$, obtain by Tonelli-Fubini's theorem, that
\begin{eqnarray*}
\er[h(Z)]&=&h(0)+\int_0^{\infty} h'(u) \; \pr(Z>u)  \;du \nonumber\\
&=& h(0)+\int_0^{x_0} h'(u) \; \pr(Z>u)  \;du + \int_{[x_0,\infty)}\int_{x_0}^x  \pr(Z>u)\;du \;\rho_h(dx),\\
&=& h(0)+\int_0^{x_0} h'(y)  \;\pr(Z>u) \;du -   h'(x_0) \int_{0}^{x_0} \pr(Z>u) \;du
+\int_{[x_0,\infty)} \int_{0}^{x} \pr(Z>u)\;du  \;\rho_h(dx).
\label{reph} \end{eqnarray*}
Since the first three terms in last equality are finite, deduce that
\begin{equation}\label{hh1}
\er[h(Z)]<\infty \Longleftrightarrow \int_{x_0}^{\infty} \int_0^{x} \pr(Z>u)\;du \;\,\rho_h(dx) <\infty.
\end{equation}
Now, observe that for all $\lambda \geq 0$,
\begin{equation}\label{contor}
(1-e^{-1})\,    (\lambda\wedge 1)  \leq \left(1-e^{-\lambda}\right)\leq (\lambda\wedge 1).
\end{equation}
Thus, write for all $x>0$,
\begin{equation}\label{contor1}
(1-e^{-1})\,x\,(\phi(1/x)\wedge 1)  \leq x\,\left(1-e^{-\phi(1/x)}\right) = \int_0^\infty e^{- u/x}\;\pr(Z>u)\;du \leq  x\;(\phi(1/x)\wedge 1),
\end{equation}
\begin{equation}\label{hh2}
\int_0^\infty e^{-u/x}\;\pr(Z>u)\;du \geq  \int_0^{x} e^{-u/x}\pr(Z>u)\;du \geq  e^{-1}\,\int_0^x \pr(Z>u) \;du,
\end{equation}
and deduce that
\begin{eqnarray}
\int_0^\infty e^{-u/x}\;\pr(Z>u)\;du&=&\int_0^{x}e^{-u/x}\;\pr(Z>u)\;du+\int_{x}^\infty e^{-u/x}\;\pr(Z>u)\;du \nonumber\\
 &\leq& \int_0^x \pr(Z>u) \;du + \pr(Z>x)\,\int_{x}^\infty e^{-u/x}\; du= \int_0^x \pr(Z>u)\; du+ e^{-1}\;x\;\pr(Z>x)\nonumber\\
&\leq& \left(1+e^{-1}\right) \int_0^x \pr(Z>u)\; du. \label{hh3}
\end{eqnarray}
From \eqref{contor1},  \eqref{hh2}  \eqref{hh3},  deduce that
$$\frac{e-1 }{e+1}\;(\phi(1/x)\wedge 1) \geq \int_0^x \pr(Z>u)\; du \leq e \;(\phi(1/x)\wedge 1),$$
and from equivalence  \eqref{hh1},  obtain the one in \eqref{ffin3}.\\

\noindent 2) Let $x_1:=\inf\{x>x_0\;;\; \phi(1/x)>1\}>0$. By \eqref{ffin3} and by the fact that
$$x\mapsto x\,\phi(1/x)=\dd +\int_0^\infty e^{-u/x} \Pi(u,\infty)\;du\;\; \mbox{is non-decreasing and starts from $\lim_{x\to 0+}x\;\phi(1/x)=\dd$}, $$
deduce the following steps:\\

a)  Assume $x_1=\infty$, then
$$\er[h(Z)]<\infty \Longleftrightarrow \int_{[x_0,\infty)}   \phi(1/x)\,x\;\rho_h(dx)   <\infty, $$
and there is no problem of integrability of the last expression at 0 if $x_0=0$.\\

b)  Assume $x_1<\infty$,  then,
$$\er[h(Z)]<\infty \Longleftrightarrow \int_{[x_0,x_1)} \,x\;  \rho_h(dx)+ \int_{[x_1,\infty)}\phi(1/x)\,x\;  \rho_h(dx)  <\infty.$$

c) Deduce that in all cases,
\begin{equation}\label{hh4}
\er[h(Z)]<\infty \Longleftrightarrow I(x_0):=\int_{[x_0,\infty)}  \phi(1/x)\;x\;  \rho_h(dx)  <\infty.
\end{equation}

d) By representation \eqref{lev1} of $\phi$ and  Tonelli-Fubini's theorem, deduce the representation,
$$I(x_0) = \dd\, h'(x_0) + \int_{\oi} \int_{[x_0,\infty)}  (1-e^{-u/x}) \;x\;  \rho_h(dx) \;\;\Pi(du),$$
and  from \eqref{contor}, deduce that
$$I(x_0)<\infty \Longleftrightarrow \int_{\oi} \chi_h(u)  \;\Pi(du),\quad \mbox{where $\chi_h(u)$ is given by \eqref{chi}.}$$
If $x_0>0$, then as a L\'evy measure, $\Pi$ always integrates $\chi_h$ on $(0,x_0)$ and the constants on $[x_0,\infty)$. The latter gives the equivalence $(i)\Longleftrightarrow (ii)$. If $x_0=0$, $\chi_h(u)= h(u)-h(0)$. Finally, from \eqref{hh4} deduce the equivalence $(ii)\Longleftrightarrow (iii)$.
\end{proofof}
\subsection{Proofs of the main results}
\begin{proofof}{\bf Theorem \ref{gamXY}.} Using representations \eqref{logammat}  for the Gamma function and \eqref{nota} for Bernstein functions, then performing an obvious change of variable, using 1) in Lemma \ref{hats}, and the fact that   $\int_{0}^{\infty} g_{\alpha ,t}(x)dx= -\dd_\alpha$ provided by Proposition \ref{hats} 1)(i), obtain the following representation  for $\alpha \in (0,1)$ and $t>0$:
\begin{equation}\label{GATT}
G_{\alpha,t}(\lambda)= \exp\left\{\dd_\alpha\lambda+ \int_0^\infty (e^{-\lambda u}-1 +\lambda u)\, \frac{g_{\alpha ,t}(u)}{u} \,du \right\}=\exp- \left\{ \int_0^\infty (1- e^{-\lambda u})  \, \frac{g_{\alpha ,t}(u)}{u} \,du \right\},
\end{equation}
where $\dd_\alpha$ is given by \eqref{nota}. Thus, since $G_{\alpha,t}(0)=1$ and using \eqref{lev1},  retrieve the equivalences in case $\alpha \in (0,1)$:
\begin{eqnarray*}
\lambda \mapsto \left(G_{\alpha,t}\right)^r \in \CM ,\quad \forall r>0&\Longleftrightarrow&  \lambda\mapsto \int_0^\infty (1- e^{-\lambda u})  \, \frac{g_{\alpha ,t}(u)}{u} \,du \in \BF \Longleftrightarrow g_{\alpha ,t} \geq 0 \\
&\Longleftrightarrow&  G_{\alpha,t}(\lambda)= \er[e^{-\lambda X_{\alpha,t}}], \; \lambda \geq 0, \quad \mbox{and}\; X\sim \ID(\RP).
\end{eqnarray*}
By point 2) in Proposition \ref{hats}, we deduce that   $g_{\alpha ,t} \geq 0 \Longleftrightarrow t\geq 1/2$. The rest of the statements are obtained by  the reflexive relation \eqref{nice}.

\end{proofof}
\medskip

\begin{proofof}{\bf Theorem \ref{gself}.} By \eqref{brsd}, the functions $g_{\alpha ,t}(x)/x$ is the density of the L\'evy measure of $X_{\alpha ,t}, \; \alpha \in (0,1)$. By point 1)(iii) in Proposition \ref{hats} is nondecreasing if, and only if $t\geq  t_0$. The reasoning is analog for $Y_{\alpha ,t}, \; \alpha >1$.
\end{proofof}
\medskip

\begin{proofof}{\bf Corollary \ref{pidg}.} 1), 2) and 3) are a straightforward consequence of the Mellin transform representations \eqref{rab}  and of Theorems \ref{gamXY} and \ref{gself}. \\
4) For the claim on $X_{\underline{\alpha},t}$,  we only need to use \eqref{GATT},  to observe that
\begin{equation}\label{berix}
-\log \er[e^{-\lambda X_{\underline{\alpha},t}}]= \int_0^\infty (1- e^{-\lambda u})  \, \frac{h_{\underline{\alpha} ,t}(u)}{u} \,du
\end{equation}
where $h_{\underline{\alpha} ,t}(u)$ is given by \eqref{gta} and to apply point 3) of Lemma  \ref{hatam} to finally obtain that
$$-\log \er[e^{-\lambda X_{\underline{\alpha},t}}]^{1/(1+t)}= \int_0^\infty  \frac{\lambda}{\lambda+x}  \frac{\eta(x)}{x} \,dx, \quad \mbox{with}\; \eta (x):= \frac{m_{\underline{\alpha} ,t}(x)}{(1+t)}\leq 1$$
is a Bernstein function that meats the form \eqref{bobo}. The assertion for $Y_{\underline{\alpha},t}$ is shown identically.
\end{proofof}
\medskip

\begin{proofof}{\bf Proposition \ref{maj}.} Using formula \eqref{logamma1} and performing the  change of variable $u\to u/\alpha $, we write
\begin{equation}\label{eqad}
{\Gamma(\lambda +t) \over \Gamma(t+\alpha \lambda)}= \exp\left\{(1-\alpha )\Psi(t) \lambda+\int_0^\infty (e^{-\lambda u}-1 +\lambda u)\,\frac{h_{\alpha,t}(u)}{u}\,du \right\}=\er[e^{\lambda T_{t,\alpha}}] , \quad \lambda\geq 0
\end{equation}
where $h_{\alpha,t}$ is the nonnegative function given by \eqref{hat} and    $T_{t,\alpha }\sim \ID_-(\Ren)$   is such that the  identity in law \eqref{gatg} holds. As a consequence of Proposition \ref{hats}, we obtain that selfdecomposability of $T_{t,\alpha }$  is equivalent to the decreaseness of $h_{\alpha,t}$. The latter is also equivalent to $t\geq t_1$.
\end{proofof}
\medskip

\begin{proofof}{\bf Proposition \ref{converg}.} Due to the representation \eqref{berix} for the Bernstein function of $X_{\underline{\alpha},t}$, it suffices to show that the function $h_{\underline{\alpha},t}$ given by \eqref{gta}, converges as $n\to \infty$. This not difficult to obtain, because for $t>0$,
the function $u\mapsto u \,e_t(u)= e^{-tu} u/(1-e^{-u})$ is bounded on $\oi$ by a constant, say $C_t$. Since $\sum_{k=1}^n \alpha_k=1$, we have
$$e_t(u) - h_{\underline{\alpha},t}(u)=  \frac{1}{u} \sum_{k=1}^n  \frac{u}{\alpha_k} h_t\left(\frac{u}{\alpha_k}\right)\; \alpha_k \leq \frac{C_t}{u}, \quad u>0.$$
Thus, for fixed $u,\, t>0$, the bounded and increasing sequence $e_t(u) - h_{\underline{\alpha},t}(u)$ is convergent as $n\to \infty$.
\end{proofof}

\medskip

\begin{proofof}{\bf Lemma \ref{op}.}
   1) It is enough to consider the case $c\in (0,1)$. By  Tonelli-Fubini's theorem, we get
\begin{eqnarray*}
\int_0^\infty  \frac{g(x)-g(x/c)}{x} \,dx &=& \int_0^\infty \frac{1}{x} \int_{[x, x/c]} d(-g)(y) \,dx
=\int_{(0,\infty)} \int_{cy}^y \frac{dx}{x}\;\; d(-g)(y)= \big(g(\infty)-g(0+) \big) \log c.
\end{eqnarray*}
2) Both assertions $(i)$ and $(ii)$ stem from
$$\theta_c(g)(x) =  - \int_1^{1/c} x g'(x s) ds= -\int_1^{1/c}  \Theta(g)(x s) \frac{ds}{s}  \quad \mbox{and}\quad \lim_{c\to 1-}\frac{\theta_c(g)(x)}{1-c}=-\Theta(g)(x).$$
\end{proofof}
\medskip

\begin{proofof}{\bf Corollary \ref{char}.} 1) Assuming that $X\sim \SD_0(\Ren_+)\;$ (resp. $X\sim \ID_-(\mathds{\Ren}) \cap \SD_0(\Ren)$), then, as in \eqref{fcl}, the associated Bernstein function $\phi$  represented by \eqref{brsd} (respectively, the Laplace exponent $\Psi$ represented by \eqref{repsi}), satisfies the following:  if $c_1, \, c_2, \cdots,c_m\in (0,1)$, and
$c_{m,1}=(1-c_1)\ldots (1-c_m), \;c_{m,2}=(1-c_1^2)\ldots (1-c_m^2)$, then  the representations
\begin{eqnarray*}
\theta_{1/c_1}\theta_{1/c_2}\ldots \theta_{1/c_m} (\phi)(\lambda)&=&  c_{m,1} \, \dd\,\lambda+ \int_0^\infty (1-e^{-\lambda y}) \frac{\theta_{c_1}\dots \theta_{c_m}(k)(x)}{x} dx\\
\theta_{1/c_1}\theta_{1/c_2}\ldots \theta_{1/c_m} (\Psi)(\lambda)&=&  c_{m,1} \, \aaa\,\lambda+ c_{m,2}\,\bb\,\lambda^2+\int_0^\infty (e^{-\lambda x}-1+\lambda x) \frac{\theta_{c_1}\dots \theta_{c_m}(k)(x)}{x} dx
\end{eqnarray*}
are straightforward and Lemma \ref{chark} clarifies the equivalences between $(i), \;(ii), (iii)$ and $(iv)$. Due to the integrability condition \eqref{x11} (respectively \eqref{x21}), we see that in the representation \eqref{hada} of $k$, necessarily $\cc=0$ and, after some calculus, we see that the $\mu$-measure should satisfy
\begin{equation}\label{tm}
\int_{\oi}  (x^l\wedge 1) \frac{k(x)}{x}dx =\int_{\oi} a_{n,l}(y) \frac{\mu(dy)}{y}<\infty,
\end{equation}
\begin{equation}\label{anl}
a_{n,l}(y) :=\;  \frac{1}{n!} \int_0^y (x^l \wedge 1) \log^n \left(\frac{y}{x}\right)\frac{dx}{x}
\;=\;\left\{
\begin{array}{lll}
\frac{y^l}{l^{n+1}},&\;\mbox{if}& x<1\\
&&\\
\frac{\log^{n+1}(y)}{(n+1)!}+\frac{y^l}{n!}\int_{\log y}^\infty z^n e^{-lz} dz,&\;\mbox{if}&  x\geq 1\,,
\end{array}
\right.
\end{equation}
with $l=1$ if $X \sim  \SD_n(\Ren_+)$ and $l=2$ if $X\sim \overline{\ID}_-(\mathds{\Ren})\cap \SD_n(\Ren)$. Observing that $a_{n,l}(y)\sim  \log^{n+1}(y)/(n+1)!$, as $y\to \infty$, we recover the condition \eqref{loginteg} which is satisfied if $X\sim \ID_-(\mathds{\Ren}) \cap \SD_n(\Ren)$ due to \eqref{integg} and to the fact that $\lim_{y\to +\infty}a_{n,2}(y)/y =0$. The equivalence with $(iv)$ and $(v)$ are due to \eqref{fraco} in the case that $X\sim \SD_n(\Ren_+)$. The proof is identical in the case that $X\sim \ID_-(\mathds{\Ren}) \cap \SD_n(\Ren)$. The equivalence with $(vi)$ is provided by \eqref{fraco} and \eqref{finiten}.\\

2) The assertion is immediate due to the decomposition and the approximation observed  in point 3) of Proposition \ref{vaguesd},  taking into account \eqref{tm}.\\

3) The first assertion is evident from point 1). The conditions on the support of $\mu$  are due to \eqref{x21},  \eqref{x11} and to Lemma \ref{chark}, which read as follows:
taking $l=0$ if  $X\sim \SD_\infty(\Ren_+)$ (respectively $l=1$ if  $X\sim \ID_-(\Ren) \cap\SD^\infty(\Ren)$), then
$$\int_0^\infty x^l(x\wedge 1)\frac{k(x)}{x}dx= \int_{\oi} \left(\int_0^1 x^{l-u} dx
+ \int_1^\infty \frac{1}{x^{u-l+1}}dx \right) \nu(du)  <\infty$$
if, and only if, $support(\nu)=(0,1)$ if $l=0$ (respectively $support(\nu)=(1,2)$), and then necessarily
$$\int_{(0,1)}  \frac{1}{u(1-u)} \nu(du)<\infty,\;\; \mbox{if $l=0$ (respectively}\; \int_{(1,2)}  \frac{1}{(u-1)(2-u)} \nu(du)<\infty,\;\; \mbox{if $l=0$)}.$$
\end{proofof}
\section{Some background on  infinitely divisible and selfdecomposable distributions}\label{acco}
An equivalent definition for the infinite divisibility of a probability measure $\mu$ is, that there exists a real row-wise independent random variables  $\xi_{n,k}$, $k=1,2,\ldots,k_n$, $\;n \ge 1$, $\; k_n\nearrow \infty$,  satisfying the {\it  infinitesimally condition}
$$\lim_{n\to \infty}\max_{1\le k\le k_n}P(|\xi_{n,k}|> \epsilon)\to 0, \quad \mbox{for each $\epsilon>0$,}$$
and such that we have the limit in distribution
\begin{equation}\label{xn}
\xi_{n,1}+ \xi_{n,2}+ \ldots+ \xi_{n,k_n} \limdis \mu.
\end{equation}

Conversely, any infinitely divisible distribution can be obtained in the scheme \eqref{xn}; see  Gnedenko \& Kolmogorov \cite{GK}  or Lo\`eve \cite{love} or Parthasarathy \cite{Par1967}, for multi-dimensional spaces.
Similarly, the distribution $\mu$ is selfdecomposable, if  there exists two sequences $a_n>0, \;b_n \in \Ren$ and an infinitesimal triangular array of independent  real random variables of the form $\xi_{n,k}:=a_n\; Z_k$,  such that
\begin{equation}\label{zn}
a_n\;(Z_1+ Z_2+\ldots+ Z_n) +b_n  \limdis \mu.
\end{equation}

The terminology of selfdecomposability  is due to the fact that the limiting distribution $\mu$ corresponds to an infinitely divisible r.v. $X$ that satisfies \eqref{hhhh}. The most general case of limits  \eqref{zn}, where the $Z_k$'s are Banach space-valued random variables and the $a_k$'s  are chosen from a group of bounded linear operators on a Banach space in question, was studied in Jurek
\cite{JU1}. In \cite[Theorem 3.1 and Section 4]{JU1},  taking reals as a Banach space $\mathds{B}$, the multiplication by positive scalars as linear operators and all probability measures as a set $Q$, it is shown that
$$\mu\in\SD_m(\mathds{B}),\;\;\mbox{if, and only if, \eqref{zn} holds and}\;\; Z_k\sim\SD_{m-1}(\mathds{B}),\; k\ge 1,\; m\ge 1.$$

Conversely, each selfdecomposable distribution can be obtained via the limiting scheme \eqref{zn}.
It might be worth to remembers that all selfdecomposable ate absolutely continuous with respect to Lebesgue measure, cf. \cite[Section 3.8, p.162]{JV1983}. Additionally, if  we assume that the r.v.'s $Z_1,\;Z_2,\ldots,Z_n,\ldots$ have the same distribution, then we get, at the limit, the class $\STA(\Ren)$ of \emph{stable distributions} on $\Ren$, see the monograph of Zolotarev \cite{zolo} or \eqref{stap} for stable distributions. Furthermore, if  $a_n:= (\sigma\sqrt{n})^{-1}$ where $\sigma^2$ is the variance of $Z_k$ and $b_n=-n\er[Z_1]$, then in \eqref{zn} we get  Central Limit Theorem, i.e., $\mu$ is the standard normal distribution $N(0,1)$.\\

From the above way of reasoning we have the inclusions:
$$\mbox{(normal distributions)}\; \subsetneq  \STA(\Ren) \subsetneq  \SD_0(\Ren) \subsetneq \ID(\Ren).$$
The class $\SD_0(\Ren)$ is quite large and contains among others $\chi^2$, Fisher, gamma, log-gamma, etc;  see Jurek \cite{JU10}.\\

In \cite[Corollary 2.11]{JU2}, Jurek showed that when taking his operator equal to the identity, we have $X \sim\SD_m(\Ren),\; \ m \ge 0$, if, and only if,  there exists a L\'evy process ${(Y_t)}_{t\geq 0}$, such that $\mathbb{E}[\log^{m+1}(1+|Y_1|)]<\infty$ and such that we have the integral stochastic representation
\begin{equation}\label{verv}
X \simdis \int_0^\infty e^{-t}dY_{r(t)}, \quad r(t)=\frac{t^{m+1}}{(m+1)!}.
\end{equation}
Note that these integral representations allow  descriptions of classes $\SD_m(\Ren)$ in terms of characteristic functions, see \cite[Theorem 3.1]{JU2}, as Urbanik \cite{urb1972, urbanik1} obtained by the extreme point method.  For $m=0$, the constructed  L\'evy process $Y$ in \eqref{verv} is such that
$$Y(t+s)-Y(t)= e^{-s}V_{e^{-t}}+V_{e^{-s}}, \quad s,\;t>0, \quad \mbox{where $V_c$ is given by \eqref{xcy}},$$
and the random integral characterization of selfdecomposable distributions  is from Jurek \& Vervaat \cite[pp. 252-253]{JV1983}. The process $Y$ is  coined as the \emph{background driving L\'evy process} of $X$, in short, BDLP. Other constructions of BDLP's for selfdecomposable random variables  are given in Jeanblanc, Pitman \& Yor \cite{JPY}.

\bigskip

\noindent {\bf Declarations of interest:} none.\\
\noindent {\bf Acknowledgment:} 
The authors are grateful to Cyril Banderier form  Paris-Nord, who kindly helped them with Maple numerical computations of the universal constant $t_1$ in \eqref{t1}. Warm thanks to the Polish team of the department of mathematics of Wroc$\l$aw for their welcome to the first author during his several visits, they were a major source of inspiration.  \\

\end{document}